\documentclass[11pt]{article}
\usepackage{amssymb}
\usepackage{amsmath}
\usepackage{amsthm}
\usepackage[colorlinks,
linkcolor=blue,
anchorcolor=blue,
citecolor=red]{hyperref}
\usepackage{amsfonts}
\usepackage{color}
\usepackage{epsfig}
\usepackage{graphics}                 % Packages to allow inclusion of graphics
\usepackage{graphicx}
\usepackage{float}
\usepackage{epsf} \usepackage{epstopdf}
\usepackage{mathrsfs}
\usepackage{subfigure}
\usepackage{cite}

\usepackage{xcolor}
\usepackage{ulem}
\definecolor{myGreen}{rgb}{0.9, 0.99, 0.9}%{0.76,0.93,0.63}
\newtheoremstyle{thmm}{1.5ex plus 1ex minus .2ex}{1.5ex plus 1ex minus
.2ex}{\it}{}{\bfseries}{}{1em}{}
\theoremstyle{thmm}
\newtheorem{theorem}{Theorem}[section]
\newtheorem{lemma}{Lemma}[section]

\newtheorem{remark}{Remark}[section]
\newtheorem{algorithm}{Algorithm}[section] 

\allowdisplaybreaks \textwidth  6.0in \textheight 9in
\newcommand{\nn}{\nonumber}
\def\A{\mathbf{A}}

\def\n{\mathbf{n}}
\def\curl{\mathbf{curl}}
\begin{document}

\date{\today}
\allowdisplaybreaks{}
\title{\bf A generalized scalar auxiliary variable method for the time-dependent Ginzburg-Landau equations\thanks{This work was completed during the author's visit to Peking University and supported by National Natural Science Foundation of China (Grant No. 12126318 \& 12126302).}}
\author{Zhiyong Si\footnote{School of Mathematics and Information Science,
		Henan Polytechnic University, 454003, Jiaozuo, P.R. China. {\tt sizhiyong@hpu.edu.cn} (Z. Si).}
}

\maketitle

\begin{abstract}
	This paper develops a generalized scalar auxiliary variable (SAV) method for the time-dependent Ginzburg-Landau equations. The backward Euler is used for discretizing the temporal derivative of the time-dependent Ginzburg-Landau equations. In this method, the system is decoupled and linearized to avoid solving the non-linear equation at each step. The theoretical analysis proves that the generalized SAV method can preserve the maximum bound principle and energy stability, which is confirmed by the numerical results. It shows that the numerical algorithm is stable.
	\maketitle
\end{abstract}

\section{Introduction}
\setcounter{equation}{0}
The phenomenological Ginzburg-Landau (GL) complex superconductivity model is given by Ginzburg and Landau in 1950s, which describes the phenomenon of vortex structure in the superconducting/normal phase transitions. The GL equations are important model in superconducting theory. The time-dependent  GL model was derived by Gor'kov and \'{E}liashberg in \cite{GE}. Due to the highly non-linear nature of the GL model, complex energy landscape and the singular dynamic response of its solution to external conditions, the analysis and research are limited to the special cases. Numerical simulation is an important tool to study GL model, which provides further theoretical insights for superconducting phenomena. The time-dependent GL model is given by the differential equations as follows \cite{GE}
\begin{align}\label{equa}
	\left\{\begin {array}{lll}
	 \eta \psi_t+i\eta \Phi\psi +(\frac{i}{\kappa}\nabla+\A)^2\psi+|\psi|^2\psi-\psi=0,                                                          \\
	  \frac{\partial \A}{\partial t}+\curl\ \curl\  \A-\nabla(\nabla\cdot \A)+\frac{i}{2\kappa}(\psi^*\nabla\psi-\psi\nabla\psi^*)+|\psi|^2\A=\curl\ H. 
	\end{array}\right.& x\in \Omega,
\end{align}
The boundary conditions are given by
\begin{align*}
	(\nabla \psi +\A \psi)\cdot \n=0, \curl\ \A \times \n=H\times \n,\ \A\cdot\n =0, &\mbox{ on } \Gamma.
\end{align*}
The initial conditions are given as follows
\begin{align*}
	\psi(x,0)=\psi_0(x), \A(x,0)=\A_0(x), \mbox{ in } \Omega,
\end{align*}
where $\Omega$ is a bounded domain in $\mathbb{R}^d$, $d=2,3$, $\n$ is the unit outer normal vector to $\Gamma\equiv\partial \Omega$, $\psi(t,x)$ is a complex valued function and usually referred to as the order parameter so that $|\psi(t,x)|$ gives the relative density of the superconducting electron pairs, and the normal and the pure superconducting states are characterized accordingly by $|\psi(t,x)|=0$ and $|\psi(t,x)|=1$ representing the superconducting state and the normal state, respectively, while $0<|\psi(t,x)|<1$ representing the mixed state. 
$\A(t,x)$ is a real vector potential for the total magnetic field. $H(x)$  is the applied magnetic field, viewed as a vector, $\kappa$ is the GL parameter, $\eta$ ($\eta\approx 0.8$) are positive material constants \cite{LMG}.

It should be remarked that the GL equations are gauge invariant. It is clear that a suitable gauge choice must be made first. There are several well-known gauge choices, e.g., the Coulomb gauge, the Lorentz gauge and the zero electric potential gauge.  The existence and uniqueness of solutions of the time-dependent GL equations have been considered in \cite{DU94,CHL93}. Numerical methods for the GL equation have been studied extensively. In \cite{DU94b}, Du gave a finite element method for the time-dependent GL equations under the zero electric gauge. A weakly nonlinear semi-implicit Euler Lagrange finite element method (FEM) for the two-dimensional GL equations was proposed by Chen and Hoffmann \cite{CH95}, where a suboptimal $L^2$-error estimate was presented. In \cite{MH98}, Mu and Huang presented an alternating Crank-Nicolson method for the time-dependent GL equations that leaded to two decoupled algebraic subsystems, one linear and the other one semi-linear.  Chen and Dai \cite{CD01} derived a posteriori estimate for the time-dependent GL model which provided the necessary information to modify the mesh and time-step according to the varying external magnetic field and corresponding motion of vortices. In \cite{GLW14}, Gao et al. gave a Crank-Nicolson finite element method the time-dependent GL equation, and the unconditional optimal error estimation was given also. Gao and Sun \cite{GSW15,GS18} presented linearized backward Euler Galerkin-mixed finite element method is investigated for the time-dependent GL equations under the Lorentz gauge based on introducing the induced magnetic field $\sigma=\curl \A$ as a new variable. In \cite{GS16}, Gao and Sun derived a new numerical approach to the time-dependent GL equations under the zero electric potential gauge. In \cite{LZ}, Li and Zhang proposed a decoupled and linearized FEM to solve the reformulated GL equations and presented error estimates in non-smooth domains. In \cite{GT18}, Ganesh and Thompson developed a time-space fully discrete implicit SMFEM algorithm for efficiently simulating the GL system modeling superconductivity on a class of superconducting surfaces. Wu and Sun \cite{WS18} presented the analysis of linearized Galerkin FEMs for a mixed formulation of the time-dependent GL equations under the temporal gauge. Li et al. \cite{LWZ20,LZ15} derived a Hodge decomposition method for dynamic GL equations. In \cite{DU94}, the authors presented the maximum bound principle of the time-dependent GL equation. The energy stability for the time-dependent GL equation was proved in \cite{DU94,LZ}. The maximum bound principle and energy stability are important properties for the time-dependent GL equations. But there is seldom numerical algorithm that can keep both the maximum bound principle and the energy stability. In \cite{GJX19}, Gao et al. proposed a stabilized semi-implicit Euler gauge-invariant method for the numerical solution of the time-dependent GL equations based on the finite difference method, which can preserve the point-wise boundedness of the solution and energy-stable. 
 
 It is well known that the energy stability and the maximum bound principle are very important in the phase field equations, e.g., Cahn-Hilliard equation, and Allen-Cahn equation. 
The scalar auxiliary variable (SAV) method is a popular method for keeping the energy stable \cite{SXY}. In this method, the gradient flow model is rewritten as an equivalent form, then some linearized and energy stability schemes can be developed by approximating the reformulated system. However, the discrete value of the SAV is not directly linked to the free energy of the dissipative system and may lead to inaccurate solutions if the time step is not sufficiently small.
  The relaxed SAV method for gradient flows was proposed by Jiang et al. in \cite{JZZ22}.
	Recently, Zhang and Shen \cite{ZS22} proposed a generalized SAV approach with relaxation (R-GSAV) for general dissipative systems. Many authors considered numerical algorithms for preserving the maximum bound principle in the phase field model. In \cite{SZ22}, Shen and Zhang derived the spectral element method for a generalized Allen-Cahn equation coupled with passive convection for a given incompressible velocity field, which can preserve the maximum bound principle. 
    In \cite{JLQ22}, Ju et al. developed first- and second-order linear finite difference schemes for a class of Allen-Cahn type gradient flow by combining the generalized SAV approach and the exponential time integrator with a stabilization term. Some other numerical analyses can see \cite{JLQ22b,SXY19,TanTang2022,TQ20} and the references therein. The aim of this paper is to design and analyze a generalized SAV algorithm for the time-dependent GL equations, which can preserve both the energy stability and maximum bound principle. This method is a linearized method to avoid solving the non-linear equation. The theoretical analysis shows that the generalized SAV method can preserve the maximum bound principle and energy stability of the time-dependent GL equation. The numerical results conform to our theoretical results.

\section{Functional settings}

This section presents some notations and the functional settings. For any two complex functions $u,v\in \mathcal{L}^2(\Omega)$, we denote the $\mathcal{L}^2(\Omega)$ inner product and norm as follows
\begin{align*}
	(u,v)=\int_{\Omega} u(x)v^*(x)d\Omega, \|u\|_0^2=[\mathcal{R}(u,u^*)]^{1/2},
\end{align*}
where $v^*(x)$ denotes the conjugate of the complex function $v(x)$, $\mathcal{R}u$ is the real part of complex function $u$. Let $W^{k,p}(\Omega)$ be the conventional Sobolev space defined on $\Omega$, when $p=2$ we denote the Hilbert space $W^{k,2}(\Omega)=H^k(\Omega)$. We denote complex Sobolev $\mathcal{H}^k(\Omega)=\{u+iv; u,v\in H^k(\Omega)\}$ and $\mathbf{H}^k(\Omega)=[H^k(\Omega)]^d$ for the vector-valued function with $d=2,3$. To take into account the time-dependent GL equation, for any given $T>0$ and a given Hilbert space $B$, we define the following spaces as follows
\begin{align*}
	L^p(0,T; B)=\left\{f(x,t)\in B, \forall 0<t<T;\int_0^T\|f(\cdot,t)\|_B^pdt< +\infty  \right\}.
\end{align*}
We also denote
\begin{align*}
	\mathbf{V}=\mathbf{L}^\infty(0,T; \mathbf{H}^1(\Omega))\cap \mathbf{H}^1(0,T; \mathbf{L}^2(\Omega)),\\
	\mathcal{V}=\mathcal{L}^\infty(0,T; \mathcal{H}^1(\Omega))\cap \mathcal{H}^1(0,T; \mathcal{L}^2(\Omega)).
\end{align*}

The time-dependent GL model (\ref{equa}) is not
well-posed for lacking uniqueness. However, they possess a gauge invariance property, which, among
other things, implies that the physical variables of interest are indeed uniquely determined
from (\ref{equa}). Several gauge choices were thoroughly discussed in \cite{DU94}. Here, we focus our attention on the gauge that eliminates the electric potential $\Phi$, named the zero electric gauge.
This is one of the most frequently used gauge choices in numerical simulations, see, e.g., \cite{DU94,DU94b}. 
The time-dependent GL equation under the zero electric gauge is given as follows.
\begin{align}\label{equation}
	\left\{\begin {array}{lll}
	 \eta\psi_t+(\frac{i}{\kappa}\nabla+\A)^2\psi+|\psi|^2\psi-\psi=0,                                                            \\
	  \frac{\partial \A}{\partial t}+\curl\ \curl\  \A-\nabla(\nabla\cdot \A)+\frac{i}{2\kappa}(\psi^*\nabla\psi-\psi\nabla\psi^*)+|\psi|^2\A=\curl\ H,
	\end{array}\right. & x\in\Omega.
\end{align}
 The GL free energy functional is given by \cite{DU94,LZ}
  \begin{align*}
	\mathcal{G} (\psi,\A)=    \int_{\Omega}\left(\left| \frac{i}{\kappa}\nabla\psi+\A\psi\right| ^2+\frac{1}{2}(|\psi|^2-1)^2\right) d\Omega
	+\int_{\Omega}\left(|\curl \ \A- H|^2+|\nabla\cdot \A |^2\right) d\Omega.
\end{align*}

\begin{lemma}[\cite{DU94,LZ} ]
	If $|\psi_0|\leq 1$, a.e. in $\Omega$, then $|\psi(t,x)|\leq 1$ a.e. in $(0,T]\times \Omega$.
\end{lemma}
The GL free energy functional decades as follows.
\begin{theorem}[\cite{DU94} ]
 For any $t\in (0,T]$, there holds that
 \begin{align*}
	\mathcal{G} (\psi(t,x),\A(t,x))\leq \mathcal{G} (\psi(0,x),\A(0,x)).
 \end{align*}
\end{theorem}

Then, the time-dependent GL equation satisfies an energy dissipation law as follows \cite{DU94}
\begin{align*}
	\frac{d \mathcal{G}  (\psi,\A) }{dt}=-\mathcal{K}(\psi,\A),
\end{align*}
where
\begin{align*}
	\mathcal{K}(\psi,\A)= 2\int_{\Omega} (\frac{\partial\psi}{\partial t},\frac{\partial\psi}{\partial t}) 
	+(\frac{\partial\A}{\partial t},\frac{\partial\A}{\partial t}
	)d\Omega.
\end{align*}
We introduce a SAV $r(t)=\mathcal{G} (\psi(t,x),\A(t,x))$,  and $r(0)=\mathcal{G} (\psi(0,x),\A(0,x))$, the time-dependent GL equation can be rewritten with the energy law as follows
\begin{align}\label{SAV-equ}
	\left\{\begin {array}{lll}
	\eta \frac{\partial\psi}{\partial t}+(\frac{i}{\kappa}\nabla+\A)^2\psi+|\psi|^2\psi-\psi=0,                                  \\
	\frac{d r}{d t}= -\frac{r(t)}{\mathcal{G} (\psi,\A)} \mathcal{K}(\psi),    \\
	\frac{\partial \A}{\partial t}+\curl\ \curl\ \A-\nabla(\nabla\cdot \A)+\frac{i}{2\kappa}(\psi^*\nabla\psi-\psi\nabla\psi^*)+|\psi|^2\A=\curl\ H,
	\end{array}\right. x\in\Omega.
\end{align}
\begin{remark}
	There are three types of SAV methods for the phase field problem, the original SAV method \cite{SXY}, the relaxed SAV method \cite{JZZ22}, and the modified SVA method \cite{HF22,HSY}. In the original SAV method for the phase field equation, the SAV is defined as $\sqrt{E+C_0}$, where $E$ is the energy, and $C_0$ is a constant such that $E+C_0>0$.
	However, the consequence is that when the time step is not sufficiently small, the modified energy can deviate far away from the original energy, leading to inaccurate solutions. The relaxed SAV approach is based on the original SAV approach which has two limitations/shortcomings: (i) it only applies to gradient flows; (ii) it requires solving two linear systems at each time step. The generalized SAV overcomes the above limitations/shortcomings while keeping the essential advantages of the original SAV approach.
\end{remark}

This paper needs the following Gagliardo-Nirenberg-Sobolev inequality \cite{GLW14} as follows
\begin{align*}
	\|u\|_{L^p}\leq C\|u\|_{H^1}, \mbox{ for } 1\leq p\leq 6.
\end{align*}

\begin{lemma}[Discrete Gronwall's lemma \cite{He90}] 
	Let $\tau $, $B$ and $a_k$, $b_k$, $c_k$, $\gamma_k$ for all integers $k>0$, be non-negative numbers such that
	\begin{align*}
		a_n+\tau\sum_{k=0}^N b_k\leq \tau\sum_{k=0}^N\gamma_ka_k+\tau\sum_{k=0}^nc_k + B, \mbox{ for all } N>0,
	\end{align*}
	suppose that $\tau \gamma_k<1$ for all $k>0$, and set $\sigma_k=(1-\tau \gamma_k)^{-1}$, then there holds that
	\begin{align*}
		a_n+\tau\sum_{k=0}^Nb_k\leq \exp\left(\tau \sum_{k=0}^N \gamma_k\sigma_k\right)\left(\tau\sum_{k=0}^N c_k+B\right).
	\end{align*}
	
\end{lemma}

\section{The first order backward Euler generalized SAV method}\setcounter{equation}{0}

This section presents the backward Euler generalized SAV method for the time-dependent GL equations. Let $0=t_0<t_1<\ldots<t_M=T$
be a uniform partition of the time interval $[0,T]$
with $t_n=n\tau$, and $M$ being a positive integer. 
Define $\zeta(t)=\frac{r(t)}{\mathcal{G} (\psi,\A)}$ and $\xi(t)= 1-(1-\zeta(t))^2$, we have 
\begin{align}\label{Seq}
	\eta \frac{\partial\psi}{\partial t}+(\frac{i}{\kappa}\nabla+\A)^2\psi+|\psi|^2\psi-\xi\psi=0,  x\in\Omega.
\end{align}
The temporal derivative $\frac{\partial \psi}{\partial t}$ will be discretized as 
\begin{align*}
	\frac{\partial \psi}{\partial t} \thickapprox \frac{\psi^n-\xi^n \psi^{n-1}}{\tau},
\end{align*}
where $\psi^n$ is an approximation of $\psi(t_n)$, $\xi^n$ is an approximation of $\xi(t_n)$ in order $\mathcal{O}(\tau^2)$, which will be given below. Then (\ref{Seq}) can be discretized as follows 
\begin{align}
	\eta \frac{\psi^n-\xi^n \psi^{n-1}}{\tau} & +(\frac{i}{\kappa}\nabla+\A^{n-1})^2\psi^n + |\psi^{n-1}|^2\psi^{n} -\xi^n \psi^{n-1} =0.
\end{align}
 If we define $\bar{\psi}^n=\frac{\psi^n}{\xi^n}$, so we can get the generalized SAV method for the time-dependent GL equations as follows.
\begin{algorithm}\label{Al3.1}
	{\bf Step 1.} Set $r^0=r(0)$ and $\psi^0= \psi(0)$, find $\bar{\psi}^n$ and $\tilde{r}^n$ by
\begin{align}\label{2.4}
	\eta \frac{\bar{\psi}^n-\psi^{n-1}}{\tau} & +(\frac{i}{\kappa}\nabla+\A^{n-1})^2\bar{\psi}^n + |\psi^{n-1}|^2\bar{\psi}^{n} -\psi^{n-1} =0, \\
	\frac{{\A}^n-\A^{n-1}}{\tau} & + \curl\ \curl\  \A^n -\nabla(\nabla\cdot \A^n) +\frac{i}{2 \kappa}(\bar{\psi}^{n*}\nabla\bar{\psi}^n-\bar{\psi}^{n}\nabla\bar{\psi}^{n*})  + |\bar{\psi}^n|^2\A^n =\curl H .\label{2.5}\\
	\frac{\tilde{r}^n-r^{n-1}}{\tau} & = -\frac{\tilde{r}^{n}}{\mathcal{G}(\bar{\psi}^{n},\A^n)}\mathcal{K}_\tau (\bar{\psi}^{n},\psi^{n-1},\A^n,\A^{n-1}),\label{2.4b}
\end{align}
where $\mathcal{K}_\tau (\bar{\psi}^{n},\psi^{n-1},\A^n,\A^{n-1})=\int_{\Omega} (\frac{ \bar{\psi}^n-{\psi}^{n-1}}{\tau},\frac{ \bar{\psi}^n-{\psi}^{n-1}}{\tau}) +(\frac{ \A^n-\A^{n-1}}{\tau},\frac{ \A^n-\A^{n-1}}{\tau})
d\Omega $. 

{\bf Step 2.} Update $\psi^n$ as follows

\begin{align}
	\zeta^n=  & \min\left\{\frac{\tilde{r}^n}{\mathcal{G} (\bar{\psi}^{n},\A^n)},1+\sqrt{3}\right\},\label{3.6b} \\	
	\psi^{n}= & \xi^n\bar{\psi}^{n}, ~\xi^n=1-(1-\zeta^n)^2.\label{2.7}
\end{align}
{\bf Step 3.} Update $r^n$ via
\begin{align}
	r^n=&\alpha_0^n\tilde{r}^n+(1-\alpha_0^n) \mathcal{G} (\psi^{n},\A^n),\alpha_0^n\in \mathcal{V}\label{3.5}
\end{align}
where 
\begin{align}
	\mathcal{V}= \left\{ \alpha \in [0,1]; \frac{r^n-\tilde{r}^n}{\tau}=\right.& -\gamma^n \mathcal{K}_\tau (\psi^{n},\psi^{n-1},\A^n,\A^{n-1})\nn\\
	&\left.+\frac{\tilde{r}^n}{\mathcal{G}(\bar{\psi}^n,\A^n)}\mathcal{K}_\tau (\bar{\psi}^{n},\psi^{n-1},\A^n,\A^{n-1})\right\}.\label{3.6}
\end{align}
with $\gamma^n\geq 0 $ to be determined so that $\mathcal{V}$ is not empty.
\end{algorithm}
\begin{remark}
 We can choose $\alpha_0^n$ and $\gamma^n$ as follows

	1. If $\tilde{r}^n=\mathcal{G}(\psi^n,\A^n)$, we set $\alpha_0^n=0$ and $\gamma^n=\frac{\tilde{r}^n\mathcal{K}_\tau (\bar{\psi}^{n},\psi^{n-1},\A^n,\A^{n-1})}{\mathcal{G}^*(\bar{\psi}^n,\A^n)\mathcal{K}_\tau (\psi^{n},\psi^{n-1},\A^n,\A^{n-1})}$.

	2. If $\tilde{r}^n>\mathcal{G}(\psi^n,\A^n)$, we set $\alpha_0^n=0$ and
	 \begin{align}\label{3.9b}
		\gamma^n=\frac{\tilde{r}^n- \mathcal{G}(\psi^n,\A^n)}{\tau \mathcal{K}_\tau (\psi^{n},\psi^{n-1},\A^n,\A^{n-1})} + \frac{\tilde{r}^n\mathcal{K}_\tau (\bar{\psi}^{n},\psi^{n-1},\A^n,\A^{n-1})}{\mathcal{G}(\bar{\psi}^n,\A^n)\mathcal{K}_\tau (\psi^{n},\psi^{n-1},\A^n,\A^{n-1})}.
	 \end{align}

	 3. If $\tilde{r}^n<\mathcal{G}(\psi^n,\A^n)$ and $\tilde{r}^n-\mathcal{G}(\psi^n,\A^n)+\tau \frac{\tilde{r}^n}{\mathcal{G}(\bar{\psi}^n,\A^n)}\mathcal{K}_\tau (\psi^{n},\psi^{n-1},\A^n,\A^{n-1})\geq 0$, we set $\alpha_0^n=0$, we set $\alpha_0^n=0$ and $\gamma^n$ given by (\ref{3.9b}).

	 4. If $\tilde{r}^n<\mathcal{G}(\psi^n,\A^n)$ and $\tilde{r}^n-\mathcal{G}(\psi^n,\A^n)+\tau \frac{\tilde{r}^n}{\mathcal{G}(\bar{\psi}^n,\A^n)}\mathcal{K}_\tau (\bar{\psi}^{n},\psi^{n-1},\A^n,\A^{n-1})< 0$, we set $\alpha_0^n=1-\frac{\tau \tilde{r}^n \mathcal{K}_\tau (\bar{\psi}^{n},\psi^{n-1},\A^n,\A^{n-1})}{\mathcal{G}(\bar{\psi}^n,\A^n)(\mathcal{G}(\psi^n,\A^n)-\tilde{r}^n)}$ and $\gamma^n=0$.
     Then, $\alpha_0^n\in \mathcal{V}$ and (\ref{3.6}) holds in all cases.
\end{remark}
\begin{remark}
	It is obvious that $\zeta(t)=\frac{r(t)}{\mathcal{G} (\psi,\A)}=1$ and $\xi(t)= 1-(1-\zeta)^2=1$ in the continuous form. In the discrete form, we choose $\zeta^n$ as an approximation of $\zeta(t_n)$ and $\xi^n$ as an approximation of $\xi(t_n)$. It will be proved that $\zeta^n= 1+\mathcal{O}(\tau)$. Here, the updating step can be seen as a correction step, which means that the generalized SAV method can be seen as a correction method.
\end{remark}

\begin{remark}
	Equation (\ref{2.4}) and (\ref{2.5}) are linearized systems, which can be easily solved without solving the nonlinear equation.
\end{remark}

\subsection{The maximum bound principle and energy stable}
\begin{theorem}\label{Thm3.1}
	If $|\psi_0(x)|\leq 1$ a.e. in $\Omega$, then there holds that $|\psi^n|\leq 1$  a.e. in $\Omega$ and \begin{align}
		r^n-r^{n-1}= -\tau \gamma^n\mathcal{K}_\tau (\psi^{n},\psi^{n-1},\A^n,\A^{n-1})\leq 0.
	\end{align}
\end{theorem}
\begin{proof}
	We prove this theorem by mathematical induction. Firstly, there holds that $|\psi^0|\leq 1$ a.e. in $\Omega$.
    
	Then, we assume $|\psi^i|<1$ for all $0<i<n$.
	Plugging (\ref{3.5}) into (\ref{3.6}), there holds that
	\begin{align}
		(\tilde{r}^n-\mathcal{G}(\psi^n,\A^n))\alpha_0^n=&\tilde{r}^n-\mathcal{G}(\psi^n,\A^n)-\tau \gamma^n\mathcal{K}_\tau (\psi^{n},\psi^{n-1},\A^n,\A^{n-1})\nonumber\\
		&+\tau \frac{\tilde{r}^n}{\mathcal{G}(\bar{\psi}^n,\A^n)}\mathcal{K}_\tau (\bar{\psi}^{n},\psi^{n-1},\A^n,\A^{n-1}).
	\end{align}
	Using (\ref{2.4b}), we can deduce that
	\begin{align*}
		\left(1+\frac{\tau}{\mathcal{G}(\bar{\psi}^n,\A^n)}\mathcal{K}_\tau (\bar{\psi}^{n},\psi^{n-1},\A^n,\A^{n-1})\right)\tilde{r}^{n}=r^{n-1}.
	\end{align*}
	Then, there holds that
	\begin{align}\label{3.13b}
	\tilde{r}^n=\frac{r^{n-1}}{1+\frac{\tau}{\mathcal{G}(\bar{\psi}^n,\A^n)}\mathcal{K}_\tau (\bar{\psi}^{n},\psi^{n-1},\A^n,\A^{n-1})}\geq 0.
	\end{align}
	Then, from (\ref{3.6b}), we derive that $\zeta^n\geq 0$, from (\ref{3.5}) we derive $\gamma^n\geq 0$. Combining (\ref{2.4b}) and (\ref{3.6}), it yields that
	\begin{align}\label{3.14b}
		r^n-r^{n-1}= -\tau \gamma^n\mathcal{K}_\tau (\psi^{n},\psi^{n-1},\A^n,\A^{n-1})\leq 0.
	\end{align}
	
	By (\ref{3.6b}) and (\ref{3.13b}), it yields that $0<\zeta^n\leq 1+\sqrt{3}$, $0\geq |[1-(1-\zeta^n)^2]|\leq 1$. Then, we can deduce that
	\begin{align*}
		|\psi^n|=\left|[1-(1-\zeta^n)^2]\bar{\psi}^{n}\right|\leq |\bar{\psi}^n|.
	\end{align*}

	Testing (\ref{2.4}) with $2\tau (|\bar{\psi}^n|^2-1)_+ \bar{\psi}^{n*}$, where $(|\bar{\psi}^n|^2-1)_+=|\bar{\psi}^n|^2-1$, if $|\bar{\psi}^n|^2-1>0$ and $(|\bar{\psi}^n|^2-1)_+=0$, otherwise, and taking the real part, it follows that
	\begin{align*}
		&2\eta\int_\Omega(|\bar{\psi}^n|^2-1)_+  \mathcal{R}\{(\bar{\psi}^{n}-\psi^{n-1})\bar{\psi}^{n*}\} d \Omega \\
		& +2\tau \int_{\Omega} \mathcal{R}\{ (\frac{i}{\kappa}\nabla+\A^{n-1})\bar{\psi}^n (-\frac{i}{\kappa}\nabla+\A^{n-1})[(|\bar{\psi}^n|^2-1)_+\bar{\psi}^{n*}]  \}d \Omega \\
	 & +2\tau \int_{\Omega}(|\bar{\psi}^n|^2-1)_+ \mathcal{R}\{(|\psi^{n-1}|^2\bar{\psi}^{n}-\psi^{n-1},\bar{\psi}^{n*})\}d\Omega= 0.
	\end{align*}

	If $|\bar{\psi}^{n}|\leq 1$, a.e. in $\Omega$, we complete the proof. Otherwise, when $|\bar{\psi}^{n}|>1$, we can deduce that
\begin{align}
	& 2\tau \int_{\Omega} \mathcal{R}\{ (\frac{i}{\kappa}\nabla+\A^{n-1})\bar{\psi}^n (-\frac{i}{\kappa}\nabla+\A^{n-1})[(|\bar{\psi}^n|^2-1)_+\bar{\psi}^{n*}]  \}d \Omega \nonumber\\
	=& 2\tau \int_{\Omega} \mathcal{R}\{ (\frac{i}{\kappa}\nabla\bar{\psi}^n) \psi^{n*} (-\frac{i}{\kappa}\bar{\psi}^n\nabla \bar{\psi}^{n*}-\frac{i}{\kappa}\bar{\psi}^{n*}\nabla \bar{\psi}^n) \}d \Omega\nonumber\\
	&+2\tau \int_{\Omega} (|\bar{\psi}^n|^2-1)_+|(\frac{i}{\kappa}\nabla+\A^{n-1})\bar{\psi}^n |^2  d \Omega\nonumber\\
	=& 2\tau \int_{\Omega} \mathcal{R}\{ \frac{1}{\kappa^2}\left(|\bar{\psi}^{n}|^2 |\nabla \bar{\psi}^n|^2+(\bar{\psi}^{n*})^2 \nabla \bar{\psi}^n\cdot \nabla \bar{\psi}^n\right) \}d \Omega\nonumber\\
	&+2\tau \int_{\Omega} (|\bar{\psi}^n|^2-1)_+|(\frac{i}{\kappa}\nabla+\A^{n-1})\bar{\psi}^n |^2  d \Omega\geq 0.\label{3.9a}
\end{align}
Using Cauchy-Schwarz inequality, the formula $2(a-b,a)=\|a\|_0^2-\|b\|_0^2+\|a-b\|_0^2$ and $|\psi^{n-1}|\leq 1$ a.e. in $\Omega$, we derive that
		 \begin{align}
			&2\eta\int_\Omega(|\bar{\psi}^n|^2-1)_+ \mathcal{R}\{ (\bar{\psi}^{n}-\psi^{n-1})\bar{\psi}^{n*}\} d \Omega+2\tau \int_{\Omega}(|\bar{\psi}^n|^2-1)_+ \mathcal{R}\{(|\psi^{n-1}|^2\bar{\psi}^{n}-\psi^{n-1},\bar{\psi}^{n*})\}d\Omega \nonumber\\
			= &\eta  \int_\Omega(|\psi^n|^2-1)_+[  |\bar{\psi}^{n}|^2-|\psi^{n-1}|^2 + | \bar{\psi}^{n}-\psi^{n-1}|^2]
			 d \Omega\nonumber\\
			 &+2\tau \int_{\Omega}(|\bar{\psi}^n|^2-1)_+( |\psi^{n-1}|^2|\bar{\psi}^{n}|^2- |\psi^{n-1}|\, | \bar{\psi}^{n}|) d\Omega\nonumber\\
			\geq  &\eta  \int_\Omega(|\psi^n|^2-1)_+[  |\bar{\psi}^{n}|^2-|\psi^{n-1}|^2 + | \bar{\psi}^{n}-\psi^{n-1}|^2]
			 d \Omega\nonumber\\
			 &+\tau \int_{\Omega}(|\bar{\psi}^n|^2-1)_+( |\psi^{n-1}|^2|\bar{\psi}^{n}|^2- |\psi^{n-1}|^2) d\Omega+\tau \int_{\Omega}(|\bar{\psi}^n|^2-1)_+( |\psi^{n-1}|^2|\bar{\psi}^{n}|^2- |\bar{\psi}^{n}|^2 ) d\Omega\nonumber\\
			 \geq  &\eta  \int_\Omega(|\psi^n|^2-1)_+[  |\bar{\psi}^{n}|^2-|\psi^{n-1}|^2 + | \bar{\psi}^{n}-\psi^{n-1}|^2]
			 d \Omega\nonumber\\
			 &+\tau \int_{\Omega}(|\bar{\psi}^n|^2-1)_+( |\bar{\psi}^{n}|^2- |\psi^{n-1}|^2) d\Omega+\tau \int_{\Omega}(|\bar{\psi}^n|^2-1)_+( |\psi^{n-1}|^2- |\bar{\psi}^{n}|^2 ) d\Omega\nonumber\\
			=  & \eta  \int_\Omega(|\bar{\psi}^n|^2-1)_+ (|\bar{\psi}^{n}|^2-|\psi^{n-1}|^2) d \Omega. \label{3.9}
		 \end{align}
By (\ref{3.9}) and $|\psi^{n-1}|\leq 1 \mbox{ a.e. in } \Omega$, we deduce that
\begin{align*}
	&2\eta\int_\Omega(|\bar{\psi}^n|^2-1)_+ \mathcal{R}\{ (\bar{\psi}^{n}-(1+\tau/\eta )\psi^{n-1})\bar{\psi}^{n*}\} d \Omega+2\tau \int_{\Omega}(|\bar{\psi}^n|^2-1)_+ |\psi^{n-1}|^2|\bar{\psi}^{n}|^2 d\Omega\\
	\geq  & \eta \int_\Omega(|\bar{\psi}^n|^2-1)_+^2  d \Omega.
 \end{align*}
		 Then, we arrive at 
		 \begin{align*}
			\eta\int_\Omega(|\bar{\psi}^n|^2-1)_+^2 d \Omega\leq &-2\tau \int_{\Omega} \mathcal{R}\{ \frac{1}{\kappa^2}\left(|\bar{\psi}^{n}|^2 |\nabla \bar{\psi}^n|^2+(\bar{\psi}^{n*})^2 \nabla \bar{\psi}^n\cdot \nabla \bar{\psi}^n\right) \}d \Omega\\
			&-2\tau \int_{\Omega} (|\bar{\psi}^n|^2-1)_+|(\frac{i}{\kappa}\nabla+\A^{n-1})\bar{\psi}^n |^2  d \Omega\leq 0.
		 \end{align*}	
Then, there holds that $(|\bar{\psi}^n|-1)_+=0$, a.e. $\Omega$, which means that
	\begin{align}
	|\bar{\psi}^{n}|\leq 1, \mbox{ a.e. in }\Omega.
	\end{align}
 Therefore, we complete the proof.
\end{proof}

\begin{remark}
	Here, $r^n<r^{n-1}$ can be seen as a modified energy stability \cite{ZS22}, which is easily preserved in the generalized SAV method. Moreover, we will give the original energy stability in the following theorem.
\end{remark}

\begin{theorem}\label{Thm3.2}
	We choose $\alpha_0^n$ and $\gamma^n$ as Algorithm \ref{Al3.1}, given $r^{n-1}>0$, we have $r^n>0$, $\zeta^n>0$ and the scheme (\ref{2.4})-(\ref{3.5}) is unconditionally energy stable in the sense that
\begin{align}
	\mathcal{G} (\psi^{n},\A^{n})\leq 	\mathcal{G} (\psi^{n-1},\A^{n-1}), \mbox{in case 1-3},
\end{align}
and
\begin{align}
	\mathcal{G} (\bar{\psi}^{n},\A^{n})\leq \mathcal{G} (\psi^{n-1},\A^{n-1}), \mbox{in case 4}.	
\end{align}
\end{theorem}
\begin{proof}
	For cases 1-3, noting $\alpha_0^n=0$, we have $r^n=\mathcal{G}(\psi^n,\A^n)$. For case 4, since $\alpha_0^n=1-\frac{\tau \tilde{r}^n \mathcal{K}_\tau (\bar{\psi}^{n},\psi^{n-1},\A^n,\A^{n-1})}{\mathcal{G}(\bar{\psi}^n,\A^n)(\mathcal{G}(\psi^n,\A^n)-\tilde{r}^n)}\in [0,1]$ and $\tilde{r}^n< \mathcal{G}(\psi^n,\A^n)$, for (\ref{3.6}) we have $r^n\leq \mathcal{G}(\psi^n,\A^n)$.

	Furthermore, in case 1-3, $\alpha_0^n=0$. From (\ref{3.5}), $r^n=\mathcal{G}(\psi^n,\A^n)$. Thanks to (\ref{3.14b}), there holds that $\mathcal{G}(\psi^n,\A^n)\leq \mathcal{G}(\psi^{n-1},\A^{n-1})$. 

	For case 4, testing (\ref{2.4}) with $2 \bar{\psi}^{n*}$ and taking the real part, it follows that
	\begin{align*}
		&2\eta\int_\Omega \mathcal{R}\{(\bar{\psi}^{n}-\psi^{n-1})\bar{\psi}^{n*}\} d \Omega  +2\tau \int_{\Omega} \mathcal{R}\{ (\frac{i}{\kappa}\nabla+\A^{n-1})\bar{\psi}^n (-\frac{i}{\kappa}\nabla+\A^{n-1})\bar{\psi}^{n*}  \}d \Omega \\
		& +2\tau \int_{\Omega} \mathcal{R}\{(|\psi^{n-1}|^2\bar{\psi}^{n}-\psi^{n-1},\bar{\psi}^{n*})\}d\Omega= 0.
	\end{align*}
Using Cauchy-Schwarz, Young's inequality, the formula $2(a-b,a)=\|a\|_0^2-\|b\|_0^2+\|a-b\|_0^2$ and $|\psi^{n-1}|, |\psi^n|\leq 1$ a.e. $\Omega$, we derive that
		 \begin{align*}
			&2\eta\int_\Omega \mathcal{R}\{ (\bar{\psi}^{n}-\psi^{n-1})\bar{\psi}^{n*}\} d \Omega+2\tau \int_{\Omega}  |\psi^{n-1}|^2|\bar{\psi}^{n}|^2 -\psi^{n-1} \bar{\psi}^{n} d\Omega\\
			\geq &\eta  \int_\Omega [  |\bar{\psi}^{n}|^2-|\psi^{n-1}|^2]	 + | \bar{\psi}^{n}-\psi^{n-1}|^2		 d+2\tau \int_{\Omega} |\psi^{n-1}|^2|\bar{\psi}^{n}|^2-|\psi^{n-1}||\bar{\psi}^{n}| d\Omega\\
			\geq  &\eta  \int_\Omega (|\bar{\psi}^{n}|^2-|\psi^{n-1}|^2) d \Omega+\tau \int_{\Omega} 2 |\psi^{n-1}|^2|\bar{\psi}^{n}|^2 -|\psi^{n-1}|^2-|\bar{\psi}^{n}|^2+\eta/\tau | \bar{\psi}^{n}-\psi^{n-1}|^2 d\Omega\\
			\geq  &\eta  \int_\Omega (|\bar{\psi}^{n}|^2-|\psi^{n-1}|^2) d \Omega+\tau \int_{\Omega} 2 |\psi^{n-1}||\bar{\psi}^{n}| -|\psi^{n-1}|^2-|\bar{\psi}^{n}|^2+\eta/\tau | \bar{\psi}^{n}-\psi^{n-1}|^2 d\Omega\\
			= &\eta  \int_\Omega (|\bar{\psi}^{n}|^2-|\psi^{n-1}|^2) d \Omega+\tau \int_{\Omega} \eta/\tau | \bar{\psi}^{n}-\psi^{n-1}|^2-(| \bar{\psi}^{n}|-|\psi^{n-1}|)^2 d\Omega.
		 \end{align*}
It is obviously that
\begin{align*}
	\int_{\Omega} \mathcal{R}\{ (\frac{i}{\kappa}\nabla+\A^{n-1})\bar{\psi}^n (-\frac{i}{\kappa}\nabla+\A^{n-1})\bar{\psi}^{n*}  \}d \Omega= \|(\frac{i}{\kappa}\nabla+\A^{n-1})\bar{\psi}^n\|_0^2 \geq 0.
\end{align*}
When $\tau \leq \eta$, using Triangle inequality, we deduce that 
\begin{align*}
	\tau \int_{\Omega} \eta/\tau | \bar{\psi}^{n}-\psi^{n-1}|^2-(| \bar{\psi}^{n}|-|\psi^{n-1}|)^2 d\Omega \geq 0.
\end{align*}
It follows that
\begin{align}\label{3.12a}
	&\int_\Omega |\bar{\psi}^{n}|^2-|\psi^{n-1}|^2 d \Omega\leq 0,\\
	&\int_\Omega |\bar{\psi}^{n}|^2-1 d \Omega\leq \int_\Omega |\psi^{n-1}|^2-1 d \Omega.\label{3.13}
\end{align}

	Multiplying   (\ref{2.4}) by $2( \bar{\psi}^{n*}-\psi^{(n-1)*})$ and taking the real part, we arrive at
	\begin{align*}
		&\frac{2\eta}{\tau}\int_\Omega 	\mathcal{R}\{(\bar{\psi}^n-\psi^{n-1})( \bar{\psi}^{n*}-\psi^{(n-1)*})\} d\Omega\\
		& +2\int_{\Omega }\mathcal{R}\left\{\left((\frac{i}{\kappa}\nabla+\A^{n-1})\bar{\psi}^n,(-\frac{i}{\kappa}\nabla+\A^{n-1})(\bar{\psi}^{n*}-\psi^{(n-1)*}) \right)\right\}d\Omega \\
		 & +2\int_{\Omega}\mathcal{R}\left\{ (|\psi^{n-1}|^2\bar{\psi}^{n}-\psi^{n-1},
		\bar{\psi}^{n*}-\psi^{(n-1)*} )\right\}d\Omega =0.
	\end{align*}
	Using the formula $2(a-b,a)=\|a\|_0^2-\|b\|_0^2+\|a-b\|_0^2$, we derive that
	\begin{align*}
		  & \mathcal{R}\left\{\left((\frac{i}{\kappa}\nabla+\A^{n-1})\bar{\psi}^n,(-\frac{i}{\kappa}\nabla+\A^{n-1})(\bar{\psi}^{n*}-{\psi}^{(n-1)*}) \right)\right\} \\
		= & \left|(\frac{i}{\kappa}\nabla+\A^{n-1})\bar{\psi}^n \right|^2-\left|(\frac{i}{\kappa}\nabla+\A^{n-1}){\psi}^{n-1} \right|^2+\frac{1}{2}\left|(\frac{i}{\kappa}\nabla+\A^{n-1})(\bar{\psi}^{n}-{\psi}^{n-1}) \right|^2 \\
		= & \frac{1}{2}\left|(\frac{i}{\kappa}\nabla+\A^{n})\bar{\psi}^n \right|^2-\frac{1}{2}\left|(\frac{i}{\kappa}\nabla+\A^{n-1}){\psi}^{n-1} \right|^2+\frac{1}{2}\left|(\frac{i}{\kappa}\nabla+\A^{n-1})(\bar{\psi}^{n}-{\psi}^{n-1}) \right|^2 \\
		  & +\frac{1}{2}\mathcal{R}\left\{\left(\frac{i}{\kappa}\nabla \bar{\psi}^n, (\A^{n-1}-\A^{n}) \bar{\psi}^{n*}\right) + \left(\frac{i}{\kappa}\nabla\bar{\psi}^{n*}, (\A^{n-1}-\A^{n}) \bar{\psi}^{n}\right) \right\}                          \\
		  & +\frac{1}{2}\left(((\A^{n-1})^2\bar{\psi}^n,\bar{\psi}^{n*})-((\A^n)^2\bar{\psi}^n,\bar{\psi}^{n*})\right) \\
		= & \frac{1}{2}\left|(\frac{i}{\kappa}\nabla+\A^{n})\bar{\psi}^n \right|^2-\frac{1}{2}\left|(\frac{i}{\kappa}\nabla+\A^{n-1}){\psi}^{n-1} \right|^2+\frac{1}{2}\left|(\frac{i}{\kappa}\nabla+\A^{n-1})(\bar{\psi}^{n}-{\psi}^{n-1}) \right|^2 \\
		  & +\mathcal{R}\left\{\left(\frac{i}{\kappa}\nabla \bar{\psi}^n, (\A^{n-1}-\A^{n}) \bar{\psi}^{n}\right)  \right\}+\left\{\frac{1}{2}\left(((\A^{n-1})^2\bar{\psi}^n,\bar{\psi}^{n*})-((\A^n)^2\bar{\psi}^n,\bar{\psi}^{n*})\right)\right\}.
	\end{align*}
	Then, it yields that
	\begin{align}
		2\tau\eta\|D_\tau \bar{\psi}^n \|_0^2 & +\left\|(\frac{i}{\kappa}\nabla+\A^{n})\bar{\psi}^n \right\|_0^2-\left\|(\frac{i}{\kappa}\nabla+\A^{n-1})\psi^{n-1} \right\|_0^2+\left\|(\frac{i}{\kappa}\nabla+\A^{n-1})(\bar{\psi}^{n}-\psi^{n-1}) \right\|_0^2\nonumber \\
	 & +2\int_\Omega\mathcal{R}\left\{\left(\frac{i}{\kappa}\nabla \bar{\psi}^n, (\A^{n-1}-\A^{n}) \bar{\psi}^{n*}\right) +\left(((\A^{n-1})^2\bar{\psi}^n,\bar{\psi}^{n*})-((\A^n)^2\bar{\psi}^n,\bar{\psi}^{n*})\right) \right\}d\Omega\nonumber         \\
	 & +2\int_\Omega\mathcal{R}\left\{ |\psi^{n-1}|^2\bar{\psi}^{n}-\psi^{n-1},\bar{\psi}^{n*}-\psi^{(n-1)*} \right\}d\Omega=0.\label{2.6}
	\end{align}

	Taking the inner product of  $2(\A^n-\A^{n-1})$ with (\ref{2.5}), we arrive at
	\begin{align*}
		2\tau \|D_\tau \A^n||_0^2 & + 2(\curl\  \A^n-H, \curl\  (\A^n-\A^{n-1}))\nn\\
 & +2\mathcal{R}\left\{(\frac{i}{\kappa}\nabla\bar{\psi}^n,\bar{\psi}^{n*} (\A^n-\A^{n-1}))\right\}+2\left(|\bar{\psi}^n|^2\A^n,(\A^n-\A^{n-1})\right)=0.
	\end{align*}
	Using the formula $2(a-b,a)=\|a\|_0^2-\|b\|_0^2+\|a-b\|_0^2$, we deduce that
	\begin{align}
		2\tau \|D_\tau \A^n||_0^2 & + \|\curl\  \A^n-H\|_0^2-\|\curl\  \A^{n-1}-H\|_0^2+\|\curl\  \A^n- \curl\  \A^{n-1}\|_0^2\nonumber\\
		&+\|\nabla\cdot \A^n\|_0^2-\|\nabla \cdot\A^{n-1}\|_0^2+\|\nabla\cdot(\A^n-\A^{n-1})\|_0^2\nn \\
	 & +2\int_\Omega \mathcal{R}\left\{(\frac{i}{\kappa}\nabla\bar{\psi}^n,\bar{\psi}^{n*} (\A^n-\A^{n-1}))\right\}d\Omega+2\left(|\bar{\psi}^n|^2\A^n,(\A^n-\A^{n-1})\right)=0.
	\end{align}

	Summing of (\ref{2.7}) and (\ref{2.6}), it yields that
	\begin{align}
		2\tau\eta\|D_\tau \bar{\psi}^n \|_0^2 & +2\tau \|D_\tau \A^n||_0^2+\left\|(\frac{i}{\kappa}\nabla+\A^{n})\bar{\psi}^n \right\|_0^2-\left\|(\frac{i}{\kappa}\nabla+\A^{n-1})\psi^{n-1} \right\|_0^2\nn \\
		 & +\left\|(\frac{i}{\kappa}\nabla+\A^{n-1})(\bar{\psi}^{n}-\psi^{n-1}) \right\|_0^2+\|\curl\  \A^n-H\|_0^2-\|\curl\  \A^{n-1}-H\|_0^2\nonumber                            \\
	& +\|\curl\  \A^n- \curl\  \A^{n-1}\|_0^2 +(\| \A^{n-1} \bar{\psi}^n\|_0^2+\|\A^n\bar{\psi}^n\|_0^2) -(|\bar{\psi}^n|^2 \A^n, \A^{n-1})\nonumber                                   \\
 & +\|\nabla\cdot \A^n\|_0^2-\|\nabla \cdot\A^{n-1}\|_0^2+\|\nabla\cdot(\A^n-\A^{n-1})\|_0^2\nn \\	                            
   & + 2\int_\Omega\mathcal{R}\left\{ |\psi^{n-1}|^2\bar{\psi}^{n}-\psi^{n-1},\bar{\psi}^{n*}-\psi^{(n-1)*}\right\}d\Omega=0.\label{2.8}
	\end{align}
Using Cauchy-Schwarz inequality and Young's inequality, we have
\begin{align*}
	(\| \A^{n-1} \bar{\psi}^n\|_0^2+\|\A^n\bar{\psi}^n\|_0^2) -2(|\bar{\psi}^n|^2 \A^n, \A^{n-1})\geq 0.
\end{align*}
	Using $|\bar{\psi}^n|<1$ a.e. $\Omega$, (\ref{3.12}) and Cauchy-Schwarz inequality, there holds that
	\begin{align}\label{2.9}
		     &2 \int_\Omega\mathcal{R}\left\{ |\psi^{n-1}|^2\bar{\psi}^{n}-\psi^{n-1},\bar{\psi}^{n*}-\psi^{(n-1)*} \right\}d\Omega\nn \\
			 =& \int_\Omega   |\psi^{n-1}|^2(|\bar{\psi}^{n}|^2 -|\psi^{n-1}|^2+|\bar{\psi}^{n}-\psi^{n-1}|^2 )+|\psi^{n-1}|^2 -|\bar{\psi}^{n}|^2+|\bar{\psi}^{n}-\psi^{n-1}|^2  d\Omega\nn\\
			=  & \int_\Omega   (|\psi^{n-1}|^2-1)(|\bar{\psi}^{n}|^2-|\psi^{n-1}|^2)+(|\psi^{n-1}|^2+1)|\bar{\psi}^{n}-\psi^{n-1}|^2 d\Omega\nn\\
			\geq& \int_\Omega   (|\psi^{n-1}|^2-1)(|\bar{\psi}^{n}|^2-|\psi^{n-1}|^2) d\Omega \geq 0.
	\end{align}

	Combining (\ref{3.13}), (\ref{2.8}) and (\ref{2.9}), it follows that
	\begin{align*}
		     & \left\|(\frac{i}{\kappa}\nabla+\A^{n})\bar{\psi}^n \right\|_0^2
		+\int_{\Omega} \frac{1}{2}(|\bar{\psi}^n|^2-1)^2d\Omega+\|\curl\  \A^n-H\|_0^2+\|\nabla\cdot \A^n\|_0^2     \\
		\leq & \left\|(\frac{i}{\kappa}\nabla+\A^{n-1}) \psi^{n-1} \right\|_0^2+\int_{\Omega} \frac{1}{2}(|\psi^{n-1}|^2-1)^2d\Omega+ \|\curl\  \A^{n-1}-H\|_0^2  +\|\nabla\cdot \A^{n-1}\|_0^2.
	\end{align*}
	Therefore, the proof is completed.
\end{proof}
\subsection{Error estimation}

This subsection will give the error estimation of the numerical algorithm. We define the errors given as follows
\begin{align*}
	e_\psi^n=\psi(t_n)-\psi^n, \bar{e}_\psi^n=\psi(t_n)-\bar{\psi}^n, e_\A^n=\A(t_n)-\A^n.
\end{align*}

\begin{theorem}
	Given the initial conditions $\bar{\psi^0}=\psi^0=\psi(0), \A^0=\A(0)$, $r^0=\mathcal{G}^*(\psi^0,\A^0)$, $u(0)\in\mathcal{H}^3(\Omega)$, $\A(0)\in H^3(\Omega)$, $\psi_{tt}\in \mathcal{H}^1(\Omega)$ and $\A_{tt}\in H^1(\Omega)$, then we have, for all $0<n<M$, when $\tau$ is sufficiently small
	\begin{align*}
		\eta \|\bar{e}_\psi^n\|_0^2+\frac{\tau }{\kappa} \sum_{i=1}^n \| \nabla\bar{e}_\psi^i\|_0^2
		+\|e_\A^n\|_0^2+  \tau \sum_{i=1}^n (\|\curl\  e_\A^i\|_0^2+\|\nabla\cdot  e_\A^i\|_0^2)
		\leq & C\tau^2,\\
		\eta \|e_\psi^n\|_0^2+\frac{\tau }{\kappa} \sum_{i=1}^n \| \nabla e_\psi^i\|_0^2
		\leq & C\tau^2,
	\end{align*}
	where $C$ is a constant dependent on $T, \Omega$ but independent of $\tau$.
\end{theorem}
\begin{proof}
	Subtracting (\ref{2.4}) from (\ref{equation}) and testing it by $2\tau \bar{e}_\psi^{n*}$ and taking the real part, we deduce the error equation
	\begin{align}
		 & 2 \eta \int_\Omega \mathcal{R}\left(\bar{e}_\psi^n-e_\psi^{n-1}, \bar{e}_\psi^{n*} \right)d\Omega+2\tau \int_\Omega \mathcal{R}\left((\frac{i}{\kappa}\nabla+\A(t_n))\bar{e}_\psi^n,(-\frac{i}{\kappa}\nabla+\A(t_n))\bar{e}_\psi^{n*} \right)d\Omega \nn\\
		 &+2\tau\int_\Omega \mathcal{R} (\frac{i}{\kappa}\nabla\bar{\psi}^n,(\A(t_n)-\A(t_{n-1}))\bar{e}_\psi^{n*})d\Omega +2\tau\int_\Omega \mathcal{R} (\frac{i}{\kappa}\nabla\bar{\psi}^n,e_\A^{n-1}\bar{e}_\psi^{n*})d\Omega \nonumber\\
		 &+2\tau \int_\Omega \mathcal{R} (\A(t_n)-A(t_{n-1})\bar{\psi}^n,\frac{i}{\kappa}\nabla\bar{e}_\psi^{n*})d\Omega +2\tau \int_\Omega \mathcal{R} (\A^{n-1}\bar{\psi}^n,(\A(t_n)-\A(t_{n-1})\bar{e}_\psi^{n*})d\Omega \nn \\
		 &  +2\tau \int_\Omega \mathcal{R} (\A^{n-1}\bar{\psi}^n,e_\A^{n-1}\bar{e}_\psi^{n*})d\Omega+2\tau \int_\Omega \mathcal{R} (e_\A^{n-1}\bar{\psi}^n,\A(t_n)\bar{e}_\psi^{n*})d\Omega                                             \nn \\
		 & +2\tau\int_\Omega \mathcal{R} (\A^{n-1}\bar{\psi}^n,(\A(t_n)-\A(t_{n-1}))\bar{e}_\psi^{n*})d\Omega  +2\tau \int_\Omega \mathcal{R}(\A^{n-1}\bar{\psi}^n,e_\A^{n-1}\bar{e}_\psi^n)d\Omega  \nn  \\
		 &+\int_\Omega \mathcal{R}(|\psi(t_n))|^2\psi(t_n)-\psi(t_n),\bar{e}_\psi^{n*}) d\Omega	-2\tau\int_\Omega \mathcal{R} (|\psi^{n-1}|^2\bar{\psi}^{n}-\psi^{n-1},\bar{e}_\psi^{n*})d\Omega\nonumber\\
		 &=2\tau\int_\Omega \mathcal{R} (Tr_\psi^n,\bar{e}_\psi^{n*})d\Omega, \label{err-eq}
	\end{align}
	where $Tr_\psi^n=\frac{\partial \psi}{\partial t}|_{t_n}-D_\tau \psi^n$.

	Then, we prove the results by mathematical induction. Firstly, we can see that
	\begin{align*}
		\zeta^0=\frac{r^0}{\mathcal{G}^* (\bar{\psi}^{0},\A^0)}=1.
	\end{align*}
	It means that
	\begin{align*}
		1-\zeta^0=0.
	\end{align*}
	We assume that $|1-\zeta^k|\leq C_0 \tau$, for all $k=1,2,\ldots, n-1$, where $C_0$ is a constant dependent  on $T, \Omega$ but independent of $\tau$.

	Taking the real part of (\ref{err-eq}), it follows that
	\begin{align*}
		 &\eta \|\bar{e}_\psi^n\|_0^2-\eta \|e_\psi^{n-1}\|_0^2+ \eta\|\bar{e}_\psi^n-e_\psi^{n-1}\|_0^2+2\tau \|(\frac{i}{\kappa}\nabla+\A(t_n))\bar{e}_\psi^n\|_0^2\leq \sum_{j=1}^{9}J_i.
	\end{align*}
	Using Cauchy-Schwarz inequality, Taylor's formulation and Young's inequality, there holds that
	\begin{align*}
		J_1 & = 2\tau \left|-\mathcal{R}\{(\frac{i}{\kappa}\nabla\bar{\psi}^n,(\A(t_n)-\A(t_{n-1})) \bar{e}_\psi^n)\}\right|                                        \\
		    & \leq 2\tau \|\bar{\psi}^n\|_\infty \|\A(t_n)-\A(t_{n-1})\|_0\| \nabla \bar{e}_\psi^n \|_0\leq \frac{\tau }{8\kappa}\| \nabla \bar{e}_\psi^n \|_0^2+C\tau^3.
	\end{align*}
	By Cauchy-Schwarz inequality and Young's inequality, we derive that
	\begin{align*}
		J_2= & 2\tau \left|-\mathcal{R}\{ (\frac{i}{\kappa}\nabla\bar{\psi}^n,e_\A^{n-1}\bar{e}_\psi^n)\}\right| \leq  2\tau \|\bar{\psi}^n\|_\infty \|e_\A^{n-1}\|_0\|\nabla e_\A^{n}\|_0 \\
		\leq & \frac{\tau }{8\kappa }\| \nabla \bar{e}_\psi^n \|_0^2 +
		C\tau \|  e_\A^{n-1} \|_0^2 .
	\end{align*}
	Using Cauchy-Schwarz inequality, Taylor's formulation and Young's inequality, it yields that
	\begin{align*}
		J_3= & 2\tau \left|-\mathcal{R}\{ ((\A(t_n)-A(t_{n-1}))\bar{\psi}^n,\frac{i}{\kappa}\nabla \bar{e}_\psi^n)\}\right| \leq \frac{\tau }{8\kappa}\| \nabla \bar{e}_\psi^n \|_0^2+C\tau^3,
	\end{align*}
	and
	\begin{align*}
		J_4= & 2\tau \left|-\mathcal{R}\{(\A^{n-1}\bar{\psi}^n,(\A(t_n)-\A(t_{n-1}))\bar{e}_\psi^n)\}\right| \leq \frac{\tau }{8\kappa}\| \nabla \bar{e}_\psi^n \|_0^2+C\tau^3.
	\end{align*}
	By Cauchy-Schwarz inequality and Young's inequality, there holds that
	\begin{align*}
		J_5= & 2\tau \left|-\mathcal{R}\{ (e_\A^{n-1}\bar{\psi}^n,\A(t_n)\bar{e}_\psi^n)\}\right| \leq \frac{\tau }{8\kappa}\| \nabla \bar{e}_\psi^n \|_0^2+C\tau\| e_\A^{n-1}\|_0^2.
	\end{align*}
	Using Cauchy-Schwarz inequality, Taylor's formulation and Young's inequality, we have
	\begin{align*}
		J_6= & 2\tau \left|-\mathcal{R}\{(\A^{n-1}\bar{\psi}^n,(\A(t_n)-\A(t_{n-1}))\bar{e}_\psi^n)\}\right|\leq \frac{\tau }{8\kappa}\| \nabla \bar{e}_\psi^n \|_0^2+C\tau^3.
	\end{align*}
	By Cauchy-Schwarz inequality, Taylor's formulation and Young's inequality, we deduce that
	\begin{align*}
		J_7= & 2\tau \left|-\mathcal{R}\{(\A^{n-1}\bar{\psi}^n,e_\A^{n-1}\bar{e}_\psi^n)\}\right|                                    \\
		\leq & 2\tau \left|\mathcal{R}\{(\A(t_{n-1})\bar{\psi}^n,e_\A^{n-1}\bar{e}_\psi^n)\}\right|+
		2\tau \left|\mathcal{R}\{ (e_\A^{n-1}\psi^n,e_\A^{n-1}(\psi(t_n)-\bar{\psi}^n))\}\right|                                      \\
		\leq & C\tau \|e_\A^{n-1}\|_0^2 +\frac{\tau }{8\kappa} \|\nabla \bar{e}_\psi^n\|_0^2+\frac{\tau }{8} \|\curl e_\A^{n}\|_0^2.
	\end{align*}
	Using Cauchy-Schwarz inequality, Taylor's formulation and Young's inequality, there holds that
	\begin{align*}
		J_8= & 2\tau \left| \mathcal{R}\{ (|\psi(t_n)|^2\psi(t_n)-\psi(t_n),e_\psi^n)	-(|\psi^{n-1}|^2\bar{\psi}^{n}-\psi^{n-1},e_\psi^n)\}\right| \\
		\leq & C\tau \|e_\psi^{n-1}\|_0^2 +\frac{\tau }{8\kappa} \|\nabla \bar{e}_\psi^n\|_0^2+C\tau^3.
	\end{align*}
 	By Cauchy-Schwarz inequality, Taylor's formulation and Young's inequality, we arrive at
	\begin{align*}
		J_9=|2\tau (Tr_\psi^n,\bar{e}_\psi^n)|\leq C\tau^3+\frac{\tau}{8\kappa}\|\nabla \bar{e}_\psi^n\|_0^2.
	\end{align*}
	By (\ref{2.7}), we derive that
	\begin{align*}
		\bar{e}_\psi^{n-1}-e_\psi^{n-1}= (1-\zeta^{n-1})^2\bar{\psi}^{n-1}.
	\end{align*}
	Testing it by $2\bar{e}_\psi^{(n-1)*}$ and taking the real part, we can get
	\begin{align*}
		\|\bar{e}_\psi^{n-1}\|_0^2-\|e_\psi^{n-1}\|_0^2\leq & 2(1-\zeta^{n-1})^2\|\bar{e}_\psi^{n-1}\|_0+ \|\bar{e}_\psi^{n-1}- e_\psi^{n-1}\|_0^2  \\
		\leq                                                & 2(1-\zeta^{n-1})^2\|\bar{e}_\psi^{n-1}\|_0+ \|(1-\zeta^{n-1})^2\bar{\psi}^{n-1}\|_0^2 \\
		\leq                                                & C\tau^3 +  \tau \|\bar{e}_\psi^{n-1}\|_0^2.
	\end{align*}
	Combine the above inequalities, it follows that
	\begin{align*}
		     & \eta \|\bar{e}_\psi^n\|_0^2-\eta \|\bar{e}_\psi^{n-1}\|_0^2+ \eta \|\bar{e}_\psi^n-e_\psi^{n-1}\|_0^2+2\tau \|(\frac{i}{\kappa}\nabla+\A(t_n))\bar{e}_\psi^n\|_0^2\nn \\
		\leq &   C\tau^3 + C \tau \|\bar{e}_\psi^{n-1}\|_0^2+C\tau \|  e_\A^{n-1} \|_0^2+\frac{\tau }{8} \|\curl e_\A^{n}\|_0^2.
	\end{align*}
	Noting
	\begin{align*}
		\|(\frac{i}{\kappa}\nabla+\A(t_n))\bar{e}_\psi^n\|_0^2\geq \frac{1}{\kappa} \| \nabla\bar{e}_\psi^n\|_0^2-\|\A(t_n)\bar{e}_\psi^n\|_0^2,
	\end{align*}
	we deduce that
	\begin{align}
		     & \eta \|\bar{e}_\psi^n\|_0^2-\eta \|\bar{e}_\psi^{n-1}\|_0^2+ \eta \|\bar{e}_\psi^n-e_\psi^{n-1}\|_0^2+\frac{2\tau }{\kappa} \| \nabla\bar{e}_\psi^n\|_0^2 \nn                \\
		\leq &  C\tau^3 + C \tau \|\bar{e}_\psi^{n-1}\|_0^2+C\tau \|  e_\A^{n-1} \|_0^2+\frac{\tau }{8} \|\curl e_\A^{n}\|_0^2+C\tau \|\bar{e}_\psi^n\|_0^2.\label{3.8}
	\end{align}

	Subtracting (\ref{2.5}) from the second equation of (\ref{equation}) and testing with $2\tau e_\A^n$, there holds that
	\begin{align*}
		2\tau (D_\tau e_\A^n,e_\A^n) & + (\curl\  e_\A^n, \curl\ e_\A^n)+\mathcal{R}\left\{(\frac{i}{\kappa}\nabla\psi(t_n),\psi(t_n) e_\A^n)-(\frac{i}{\kappa}\nabla\bar{\psi}^n,\bar{\psi}^n e_\A^n)\right\}\nn \\
		                           &+(\nabla\cdot e_\A^n, \nabla\cdot e_\A^n)
		+2\tau \left(|\psi(t_n)|^2A(t_n),e_\A^n\right)-\left(|\bar{\psi}^n|^2\A^n,e_\A^n\right)=(Tr_A,e_\A^n),
	\end{align*}
	where $Tr_A=\frac{\partial A}{\partial t}|_{t_n}-D_\tau \A^n$. Using the formula $2(a-b,a)=\|a\|_0^2-\|b\|_0^2+\|a-b\|_0^2$, there holds that
	\begin{align}
		 \|e_\A^n\|_0^2&- \|e_\A^{n-1}\|_0^2+ \|e_\A^n-e_\A^{n-1}\|_0^2  + 2\tau \|\curl\  e_\A^n\|_0^2 + 2\tau \|\nabla\cdot e_\A^n\|_0^2\nn\\
		& +2\tau \mathcal{R}\left\{(\frac{i}{\kappa}\nabla\psi(t_n),\psi(t_n) e_\A^n )-(\frac{i}{\kappa}\nabla\bar{\psi}^n,\bar{\psi}^n e_\A^n)\right\}\nn \\
	 &+2\tau \left(|\psi(t_n)|^2A(t_n),e_\A^n\right)-2\tau \left(|\bar{\psi}^n|^2\A^n,e_\A^n\right)=2\tau (Tr_\A^n,e_\A^n).\label{3.10}
	\end{align}
	Using Young's and Cauchy-Schwarz's inequality, we derive that
	\begin{align*}
		     & \left|\mathcal{R}\left\{(\frac{i}{\kappa}\nabla\psi(t_n), \psi(t_n) e_\A^n)-(\frac{i}{\kappa}\nabla\bar{\psi}^n,\bar{\psi}^n e_\A^n) \right\}\right|                                                 \\
		\leq    & \left|\mathcal{R}\left\{(\frac{i}{\kappa}\nabla \bar{e}_\psi^n ,\psi(t_n) e_\A^n)\right\}\right|+ \left|\mathcal{R}\left\{(\frac{i}{\kappa}\nabla \bar{\psi}^n ,\bar{e}_\psi^n e_\A^n)\right\}\right| \\
		\leq & C\tau \|\bar{e}_\psi^{n}\|_0^2 +\frac{\tau }{8\kappa} \|\nabla \bar{e}_\psi^n\|_0^2 + \frac{\tau }{4} \|\curl e_\A^n\|_0^2.
	\end{align*}
	By Taylor's formulation, Young's and Cauchy-Schwarz's inequality, we deduce that 
	\begin{align*}
		&2\tau |\left(|\psi(t_n)|^2A(t_n),e_\A^n\right)-2\tau \left(|\bar{\psi}^n|^2\A^n,e_\A^n\right)|\\
		\leq  & 2\tau |\left(|\psi(t_n)|^2e_\A^n,e_\A^n\right)|+2\tau |\left((|\psi(t_n)|^2-|\bar{\psi}^n|^2)\A^n,e_\A^n\right)|  \\
		\leq                                                                                       & C\tau \|\bar{e}_\A^{n}\|_0^2 +\frac{\tau }{8\kappa} \|\nabla \bar{e}_\psi^n\|_0^2 + \frac{\tau }{4} \|\curl e_\A^n\|_0^2,
	\end{align*}
	and
	\begin{align*}
		2\tau |(Tr_A,e_\A^n)|\leq C\tau \|Tr_\A^n\|_0\|\curl e_\A^n\|_0\leq \frac{\tau }{4} \|\curl e_\A^n\|_0^2+C\tau^3.
	\end{align*}
	Then, there holds that
	\begin{align}
		\|e_\A^n\|_0^2-\|e_\A^{n-1}\|_0^2+\|e_\A^n-e_\A^{n-1}\|_0^2  +& \tau \|\curl\  e_\A^n\|_0^2 + \tau \|\nabla\cdot e_\A^n\|_0^2\nn \\
		\leq & C\tau \|\bar{e}_\psi^{n}\|_0^2+C\tau \|\bar{e}_\A^{n}\|_0^2 +\frac{\tau }{8\kappa} \|\nabla \bar{e}_\psi^n\|_0^2 .
		\label{3.11}
	\end{align}
	Combining (\ref{3.8}) and (\ref{3.11}), we arrive at
	\begin{align}
		     & \eta \|\bar{e}_\psi^n\|_0^2-\eta \|e_\psi^{n-1}\|_0^2+ \eta \|\bar{e}_\psi^n-e_\psi^{n-1}\|_0^2+\frac{\tau }{\kappa} \| \nabla\bar{e}_\psi^n\|_0^2\nn                                       \\
		     & +\|e_\A^n\|_0^2-\|e_\A^{n-1}\|_0^2+\|e_\A^n-e_\A^{n-1}\|_0^2 + \tau \|\curl\  e_\A^n\|_0^2+ 2\tau \|\nabla\cdot e_\A^n\|_0^2 \nn         \\
		\leq & C\tau^3 + C \tau \|\bar{e}_\psi^{n-1}\|_0^2+C\tau \|  e_\A^{n-1} \|_0^2+C\tau \|\bar{e}_\psi^{n}\|_0^2+C\tau \|\bar{e}_\A^{n}\|_0^2.\label{3.12}
	\end{align}
	Summing (\ref{3.12}) over all $n$, and using Gronwall's lemma, we have
	\begin{align}
		\|\bar{e}_\psi^n\|_0^2+\frac{\tau }{\kappa} \sum_{i=1}^n \| \nabla\bar{e}_\psi^i\|_0^2
		+\|e_\A^n\|_0^2+  \tau \sum_{i=1}^n (\|\curl\  e_\A^i\|_0^2+  \|\nabla\cdot e_\A^i\|_0^2)
		\leq C\tau^2.
	\end{align}

	At last, we will prove $|1-\zeta^n|\leq C_0 \tau$. By (\ref{2.7}), there holds that
	\begin{align*}
		|1-\zeta^n|\leq & \left| 1-\frac{r^n}{\mathcal{G}(\bar{\psi}^n,\A^n)}\right| =\left|\frac{r(t_n)}{\mathcal{G}(\psi(t_n),\A(t_n))}-\frac{r^n}{\mathcal{G}(\bar{\psi}^n,\A^n)} \right| \\
		= & \left|\frac{r(t_n)-r^n}{\mathcal{G}(\psi(t_n),\A(t_n))}-\frac{r^n[\mathcal{G}(\bar{\psi}^n,\A^n)-\mathcal{G}(\psi(t_n),\A(t_n))]}{\mathcal{G}(\bar{\psi}^n ,\A^n)\mathcal{G}(\psi(t_n),\A(t_n))}\right| \\
		\leq         & C\tau +C|\mathcal{G}(\bar{\psi}^n,\A^n)-\mathcal{G}(\psi(t_n),\A(t_n))|.
	\end{align*}
	By the definition of $	\mathcal{G}(\cdot,\cdot)$, we deduce that
	\begin{align*}
		     &\left| \mathcal{G}(\bar{\psi}^n,\A^n)-\mathcal{G}(\psi(t_n),\A(t_n)) \right|    \\
		=    &\left| \int_{\Omega } \left(\left| \frac{i}{\kappa}\nabla\bar{\psi}^n+\A^n\bar{\psi}^n\right| ^2+\frac{1}{2}(|\bar{\psi}^n|^2-1)^2\right)-\left(\left| \frac{i}{\kappa}\nabla\psi(t_n)+\A(t_n)\psi(t_n)\right| ^2+\frac{1}{2}(|\psi(t_n)|^2-1)^2\right)d\Omega \right| \\
		&+\int_{\Omega}\left(|\curl \ \A^n- H|^2+|\nabla\cdot \A^n |^2\right) d\Omega-\int_{\Omega}\left(|\curl \ \A(t_n)- H|^2+|\nabla\cdot \A(t_n) |^2\right) d\Omega\\
		\leq & C\|\nabla \bar{e}_\psi^{n}\|_0^2+C\|e_\A^n \|_0^2+ C\|\bar{e}_\psi^{n}\|_0^2\leq C\tau^2.
	\end{align*}
	Then, we can derive
	\begin{align*}
		|1-\zeta^n|\leq C_0\tau.
	\end{align*}
	Then, we finish the mathematical induction.

	Noting $\psi^{n}= [1-(1-\zeta^n)^2]\bar{\psi}^{n}$, we deduce that
\begin{align*}
	\|e_\psi^n\|_0^2 \leq & \|\psi^n-\bar{\psi}^n\|_0^2+\|\bar{e}_\psi^n\|_0^2 \leq C\tau^2+\|\bar{e}_\psi^n\|_0^2,\\
    \|\nabla e_\psi^n\|_0^2 \leq& \|\nabla (\psi^n-\bar{\psi}^n)\|_0^2+\|\nabla \bar{e}_\psi^n\|_0^2 \leq C\tau^2+\|\bar{e}_\psi^n\|_0^2.
\end{align*}
Then, there holds that
\begin{align*}
	\eta \|e_\psi^n\|_0^2+\frac{\tau }{\kappa} \sum_{i=1}^n \| e_\psi^i\|_0^2
	\leq C\tau^2.
\end{align*}
\end{proof}

\section{Numerical Results}

\begin{figure}[ht]\begin{center}
	\includegraphics[width=0.50\textwidth, height=66mm]{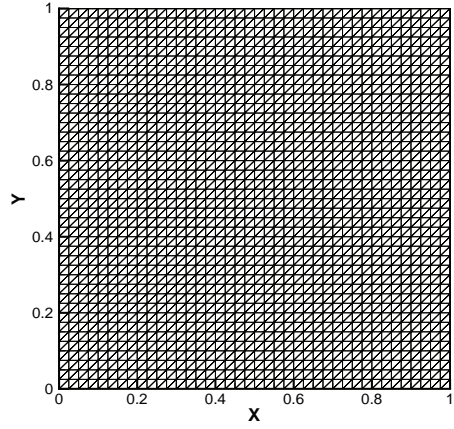}\\
	\caption{The mesh in unit square.}\label{mesh1}	
	\end{center}
\end{figure}

In this section, we present some numerical results to show the effect of the generalized SAV method for the time-dependent GL model, so the computation domain $\Omega$ is chosen as a polygon in $2$ dimension. We discretize the GL equation in space using the finite element method. 
Let $\mathcal{T}_h$ be a regular triangle partition of $\Omega$ with $\Omega=\cup_e \Omega_e $; we denote the mesh size by $h=\max_{\Omega_e\in \mathcal{T}_h}\{\mbox{Diam} \Omega_e \}$. For a given partition $\mathcal{T}_h$ , we denote $ \mathcal{V}_h^2$ and $\mathbf{V}_h^2$ as the $2$nd-order Lagrange finite element subspaces of $\mathcal{H}_1(\Omega)$ and $\mathbf{H}_n^1(\Omega)$, respectively. Here, we choose the finite element as $\mathcal{V}_h^2,~ \mathbf{V}_h^2$ for $\psi$ and $\A$. Then, we can get the finite element method of the time-dependent GL equations can be given as follows, find $(\psi_h^n, \A_h^n)\in \mathcal{V}_h^2\times \mathbf{V}_h^2$ such that
\begin{align*}
\eta \left(\frac{\bar{\psi}_h^n-\psi_h^{n-1}}{\tau} ,\tilde{\psi}_h\right)&+\left((\frac{i}{\kappa}\nabla+\A_h^{n-1})\bar{\psi}_h^n,(\frac{i}{\kappa}\nabla+\A_h^{n-1})\tilde{\psi}_h \right)\\
&+(|\psi_h^{n-1}|^2\bar{\psi}_h^n-\psi_h^{n-1},\tilde{\psi}_h)=0,& \forall \tilde{\psi}_h\in \mathcal{V}_h^r,\\
\left(\frac{\A_h^n-\A_h^{n-1}}{\tau },\tilde{\A}_h\right)& + (\curl\  \A_h^n, \curl\ \tilde{ \A}_h)+(\nabla\cdot \A_h^n,\nabla\cdot \tilde{\A}_h)+\left(|\psi_h^n|^2\A,\tilde{\A}\right)\nn\\
&+\frac{i}{2 \kappa}\left(\bar{\psi}_h^{n*}\nabla\bar{\psi}_h^n-\bar{\psi}_h^{n}\nabla\bar{\psi}_h^{n*}, \tilde{\A}_h \right)=(H,\curl\ \tilde{\A}_h), &\forall \tilde{\A}_h \in \mathbf{V}_n^r,\\
\frac{\tilde{r}_h^n-r_h^{n-1}}{\tau} & = -\frac{\tilde{r}_h^{n}}{\mathcal{G}(\bar{\psi}_h^{n},\A_h^n)}\mathcal{K}(\bar{\psi}_h^{n},\A_h^n),
\end{align*}
where $\mathcal{K}(\bar{\psi}_h^{n},\A_h^n)=\int_{\Omega} (\frac{ \bar{\psi}_h^n-\bar{\psi}_h^{n-1}}{\tau},\frac{ \bar{\psi}_h^n-\bar{\psi}_h^{n-1}}{\tau}) +(\frac{ \A_h^n-\A_h^{n-1}}{\tau},\frac{ \A_h^n-\A_h^{n-1}}{\tau})
d\Omega $.  Then, we update $\psi_h^n$ as follows
\begin{align}
	\zeta^n=  & \frac{\tilde{r}^n}{\mathcal{G} (\bar{\psi}_h^{n},\A_h^n)}, \\	
	\psi_h^{n}= & \xi^n\bar{\psi}_h^{n}, ~\xi^n=1-(1-\zeta^n)^2.
\end{align}
Then, update $r^n$ via
\begin{align}
	r^n=&\alpha_0^n\tilde{r}^n+(1-\alpha_0^n) \mathcal{G} (\psi_h^{n},\A_h^n),\alpha_0^n\in \mathcal{V}
\end{align}
where 
\begin{align}
	\mathcal{V}= \left\{ \zeta \in [0,1]; \frac{r^n-\tilde{r}^n}{\tau}=-\gamma^n \mathcal{K}(\psi_h^n,\A_h^n)+\frac{\tilde{r}^n}{\mathcal{G}(\bar{\psi}_h^n,\A_h^n)}\mathcal{K}(\bar{\psi}_h^n,\A_h^n), \gamma^n\geq 0 \right\}.
\end{align}
We choose $\alpha_0^n$ and $\gamma^n$ as follows:

	1. If $\tilde{r}^n=\mathcal{G}(\psi_h^n,\A_h^n)$, we set $\alpha_0^n=0$ and $\gamma^n=\frac{\tilde{r}^n\mathcal{K}(\bar{\psi}_h^n,\A_h^n)}{\mathcal{G}^*(\bar{\psi}_h^n,\A^n)\mathcal{K}(\psi_h^n,\A_h\mathcal{K}(\psi^n)^n)}$.

	2. If $\tilde{r}^n>\mathcal{G}(\psi^n,\A^n)$, we set $\alpha_0^n=0$ and
	 \begin{align}
		\gamma^n=\frac{\tilde{r}^n- \mathcal{G}(\psi_h^n,\A_h^n)}{\tau \mathcal{K}(\psi_h^n,\A_h^n)} + \frac{\tilde{r}^n\mathcal{K}(\bar{\psi}_h^n,\A_h^n)}{\mathcal{G}(\bar{\psi}_h^n,\A_h^n)\mathcal{K}(\psi_h^n,\A_h^n)}.
	 \end{align}

	 3. If $\tilde{r}^n<\mathcal{G}(\psi_h^n,\A_h^n)$ and $\tilde{r}^n-\mathcal{G}(\psi_h^n,\A_h^n)+\tau \frac{\tilde{r}^n}{\mathcal{G}(\bar{\psi}_h^n,\A_h^n)}\mathcal{K}(\psi_h^n,\A_h^n)\geq 0$, we set $\alpha_0^n=0$, we set $\alpha_0^n=0$ and $\gamma^n$ given by (\ref{3.9b}).

	 4. If $\tilde{r}^n<\mathcal{G}(\psi_h^n,\A_h^n)$ and $\tilde{r}^n-\mathcal{G}(\psi_h^n,\A_h^n)+\tau \frac{\tilde{r}^n}{\mathcal{G}(\bar{\psi}_h^n,\A_h^n)}\mathcal{K}(\bar{\psi}_h^n,\A_h^n)< 0$, we set $\alpha_0^n=1-\frac{\tau \tilde{r}^n \mathcal{K}(\bar{\psi}_h^n,\A_h^n)}{\mathcal{G}(\bar{\psi}_h^n,\A_h^n)(\mathcal{G}(\psi_h^n,\A_h^n)-\tilde{r}^n)}$ and $\gamma^n=0$.
	 
The code was implemented by the open source code for the finite element method, FreeFEM++ \cite{Hec}, which is a popular 2D and 3D partial differential equations (PDE) solver. It allows the authors to easily implement their own physics modules using the provided FreeFEM++ language.
	\begin{figure}[ht]
		\begin{center}
\subfigure[$\kappa=1$]{\includegraphics[width=0.4\textwidth, height=56mm]{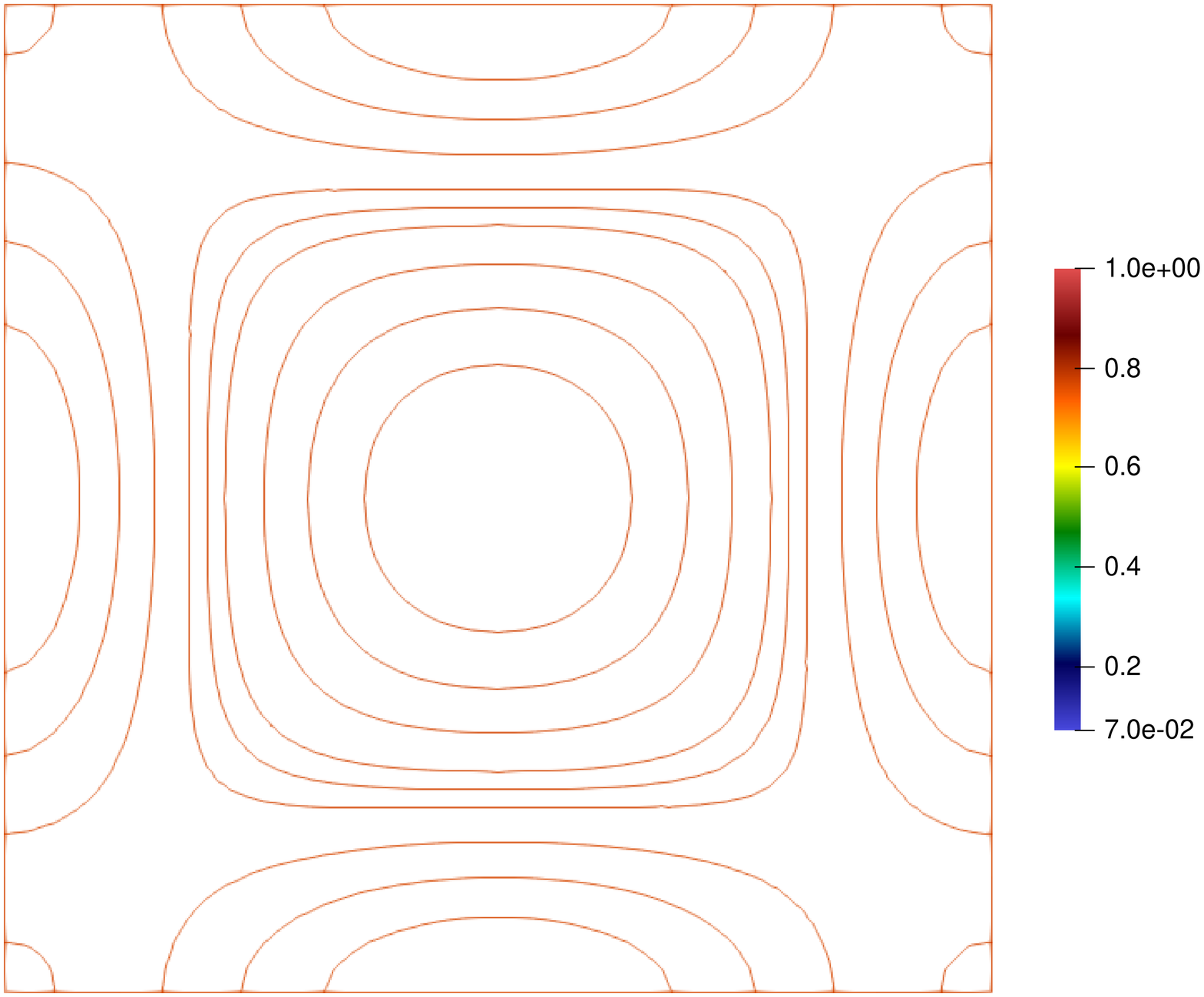}  }  \subfigure[$\kappa=10$]{\includegraphics[width=0.4\textwidth, height=56mm]{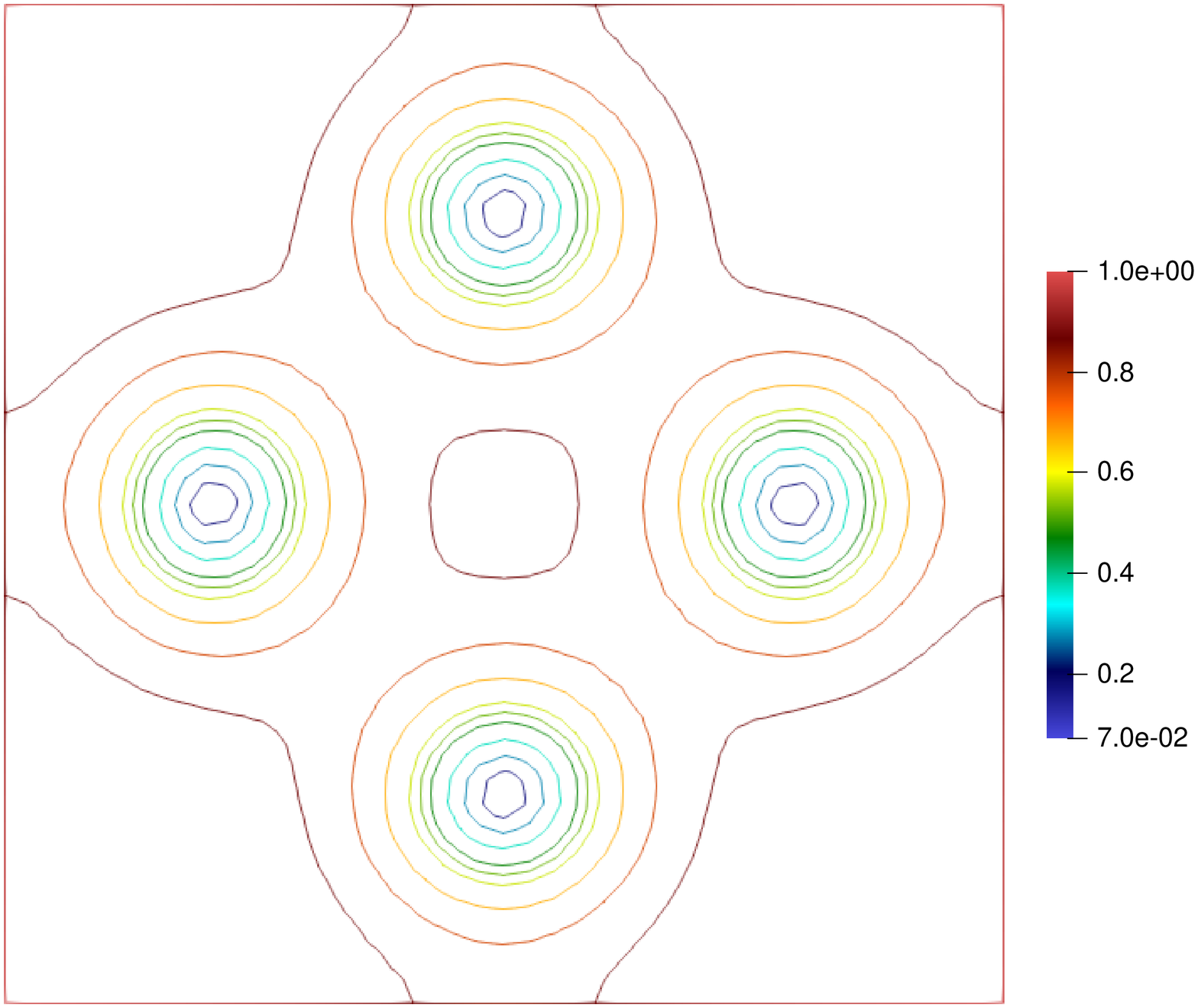}}\\
\subfigure[$\kappa=20$]{\includegraphics[width=0.4 \textwidth, height=56mm]{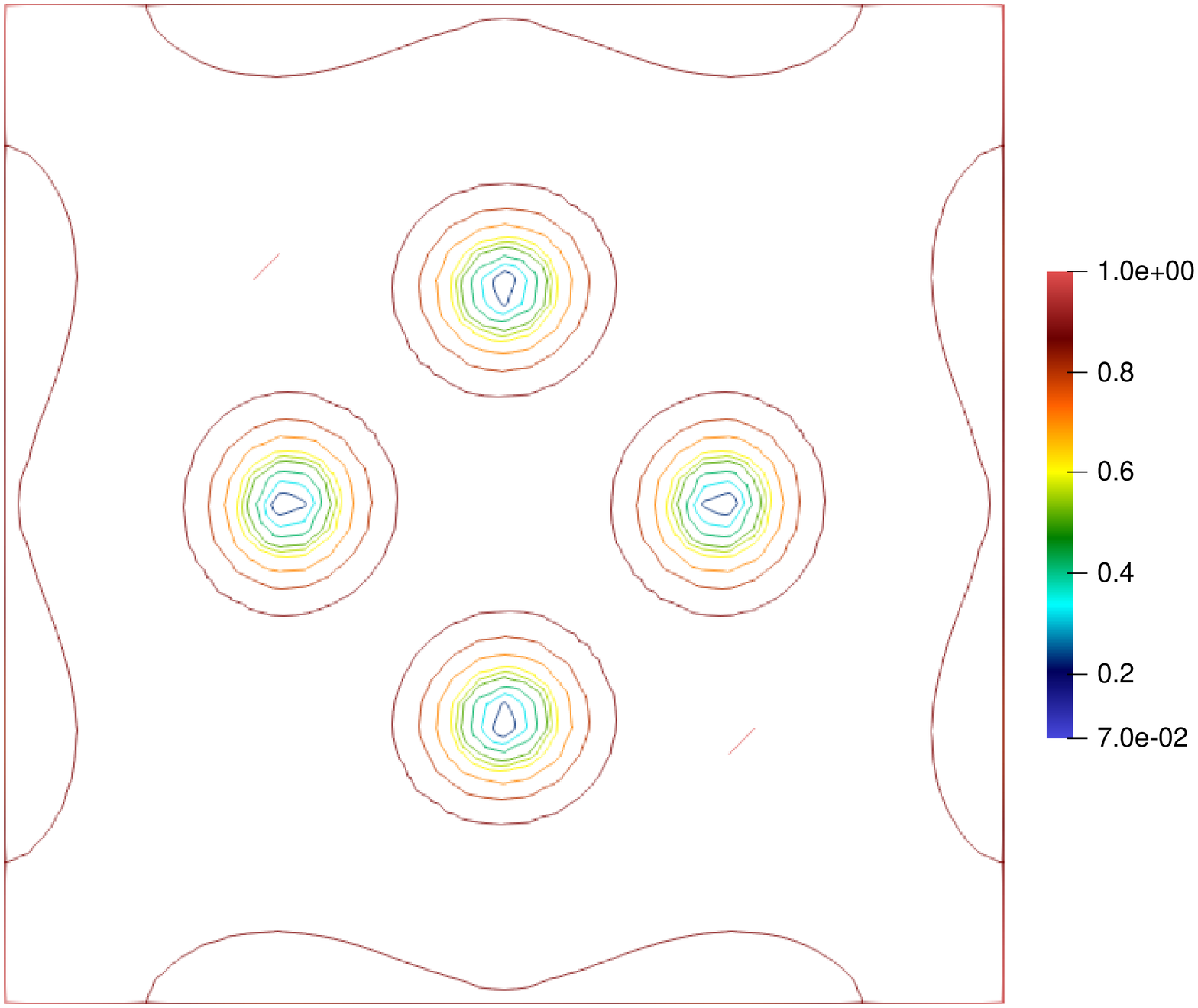}  }  \subfigure[$\kappa=50$]{\includegraphics[width=0.4 \textwidth, height=56mm]{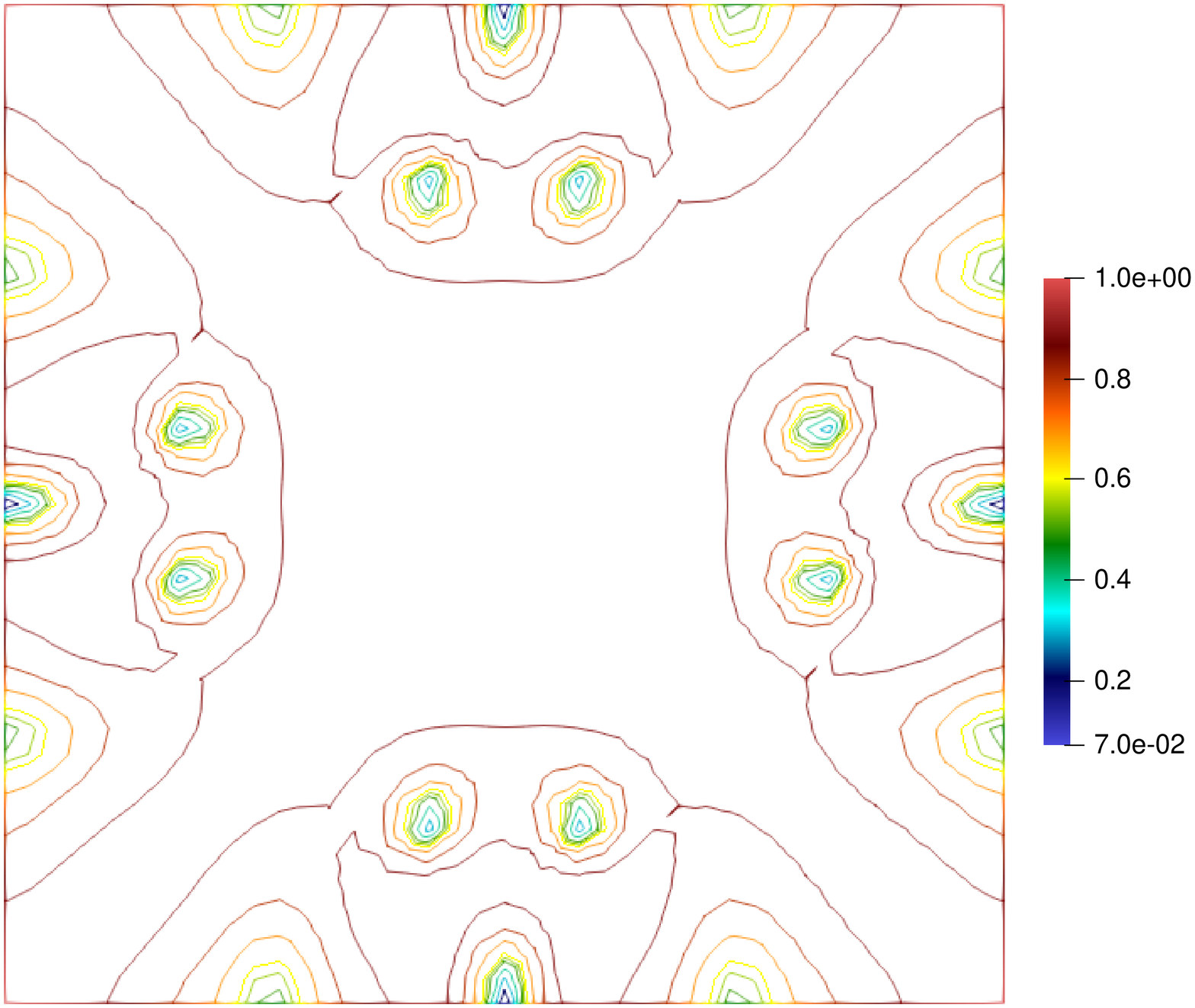}}\\
		\caption{Contour plots of $|\psi|$ in square with different $\kappa$.}\label{F1}
		\end{center}
	\end{figure}
	
\begin{remark}
	The initial condition $\phi_h^0$ should keep the maximum bound principle, i.e. $|\psi_h^0|\leq 1$ a.e. in $\Omega$. Then, we can prove that the finite element algorithm can keep the maximum bound principle also. The proof is as similar as the proof of Theorem \ref{Thm3.1}, which is omitted here. Similarly, the finite element algorithm can preserve energy stability, the proof is as similar as Theorem \ref{Thm3.2}. The numerical results show these.
 \end{remark}

\begin{remark}
	These finite element systems are linear algebraic systems. As the term $((\frac{i}{\kappa}\nabla+\A_h^{n-1})\bar{\psi}_h^n,(\frac{i}{\kappa}\nabla+\A_h^{n-1})\tilde{\psi}_h)$ and $(\curl\  \A_h^n, \curl\ \tilde{ \A}_h)+(\nabla\cdot \A_h^n,\nabla\cdot \tilde{\A}_h)$ are symmetry and positive, the linear algebraic systems are well-posed. They can be solved by the numerical algorithm for the linear algebraic system, e.g., GMRES, UMFPACK, and so on. In this paper, we use GMRES for solving finite element systems.
\end{remark}
\begin{remark}
	The maximum bound principle and the energy stable can be proved similarly as the time discrete generalized SAV method, we omit it. 
\end{remark}

\subsection{The vortex simulation in unit square}
\begin{figure}[ht]\begin{center}
			\includegraphics[width=0.66\textwidth, height=56mm]{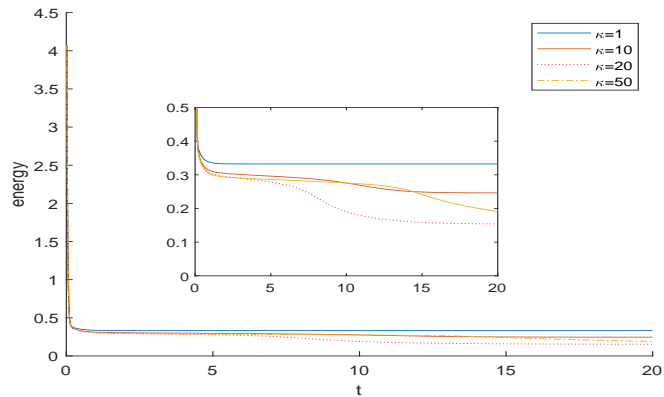}\\
			\caption{Time evolution of the discrete energy $\mathcal{G} (\psi,\A)$.}\label{F2}	
			\end{center}
		\end{figure}

\begin{figure}[ht]\begin{center}
			\includegraphics[width=0.66\textwidth, height=56mm]{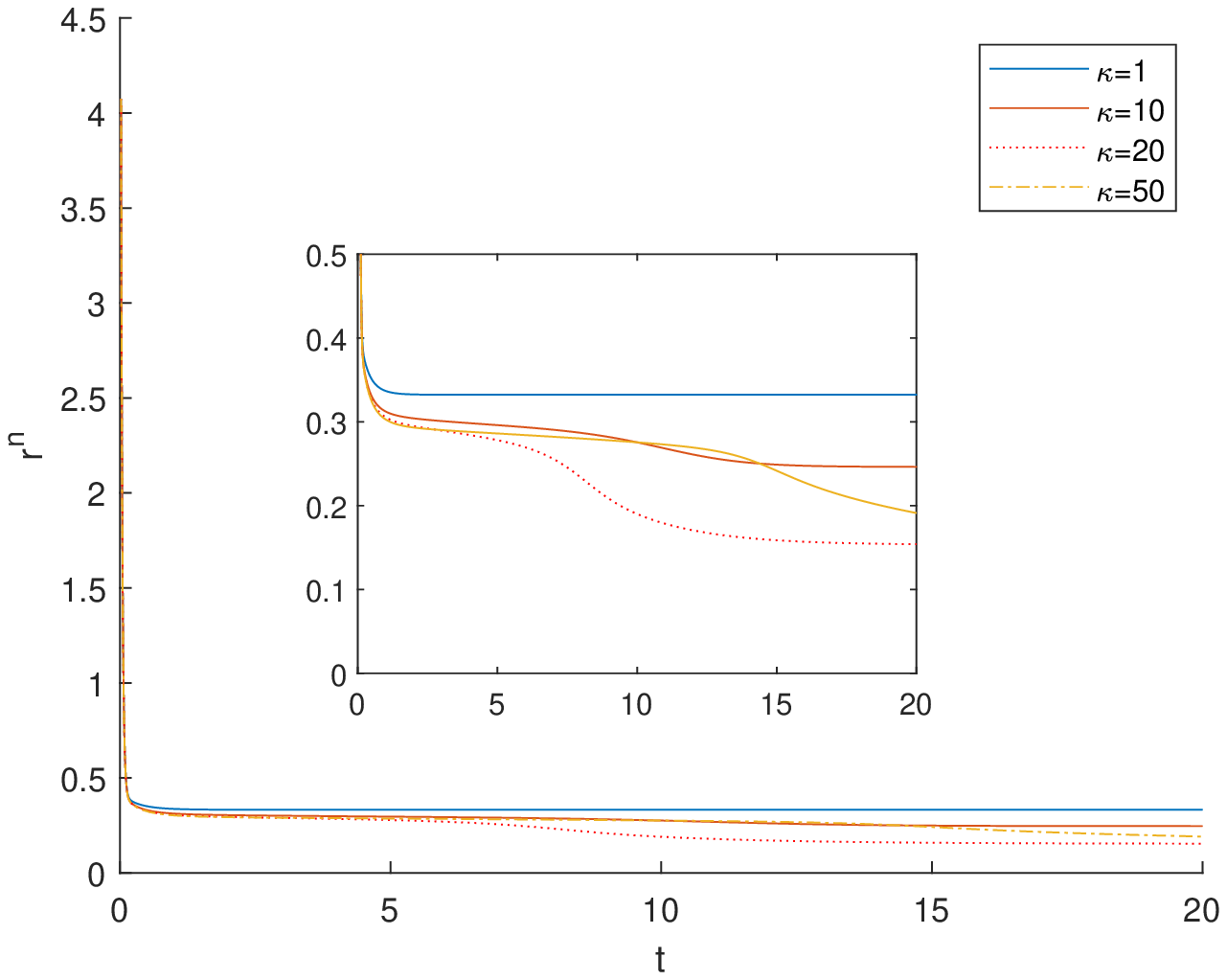}\\
			\caption{Time evolution of the SAV $r$.}\label{F3}	
			\end{center}
		\end{figure}

\begin{figure}[ht]\begin{center}
			\includegraphics[width=0.66\textwidth, height=56mm]{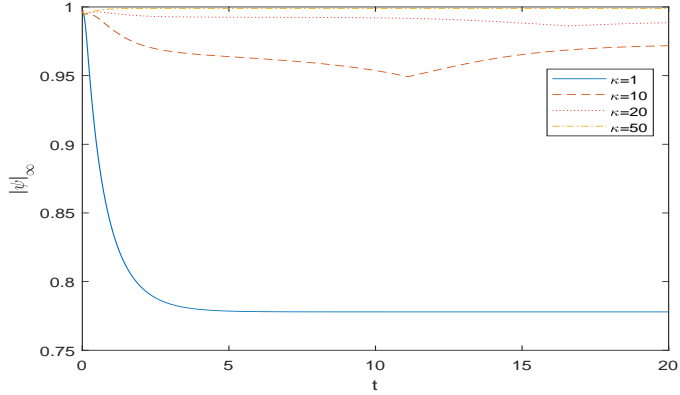}\\
			\caption{Time evolution of $|\psi|_{\infty}$.}\label{Fe13}	
			\end{center}
		\end{figure}

This subsection presents the numerical results of the vortex simulation of time-dependent GL equations
with domain $\Omega = [0, 1]\times [0, 1]$. We set the GL parameter $\kappa = 1, 10, 20$ and $50$, $\eta=1$. The initial conditions are chosen as $\psi_0 = 0.8+i0.6$ and $\A_0 = (0, 0)^T$, it means that the initial state is purely in the superconducting state. The finite mesh is chosen as the uniform triangle grid, see Figure \ref{mesh1}, and the step size is chosen as $h=\sqrt{2}/40$, $\tau=0.01$ and the final time is $T=20$.  
The applied magnetic field $H = 3.5$. Figure \ref{F1} shows the contour plots of $|\psi|$ with different $\kappa$. We can see that there is no vortex when $\kappa=1$. There are four vortexes when $\kappa=10$, which confirms the results in \cite{DU94,GLW14}. In order to show the robustness of the numerical algorithm, we show the numerical results for $\kappa=20$ and $50$. For $\kappa=20$, there are four vortexes also, but the vortexes are smaller near the center. When $\kappa=50$, there are more vortexes. Figure \ref{F2} shows the time evolution of the energy. Figure \ref{F3} shows the time evolution of $r^n$. It shows that the energy is stable conforming with the theoretical result. Figure \ref{Fe13} presents the time evolution of $|\psi|_{\infty}$, it shows that the maximum of $|\psi|$ are smaller than $1$. It confirms the theoretical analysis. The numerical results show that the generalized SAV algorithm for the time-dependent GL equation can preserve the energy stability and maximum bound principle.

\subsection{The vortex simulation in a multi-connected domain}
\begin{figure}[ht]\begin{center}
	\includegraphics[width=0.64\textwidth, height=86mm]{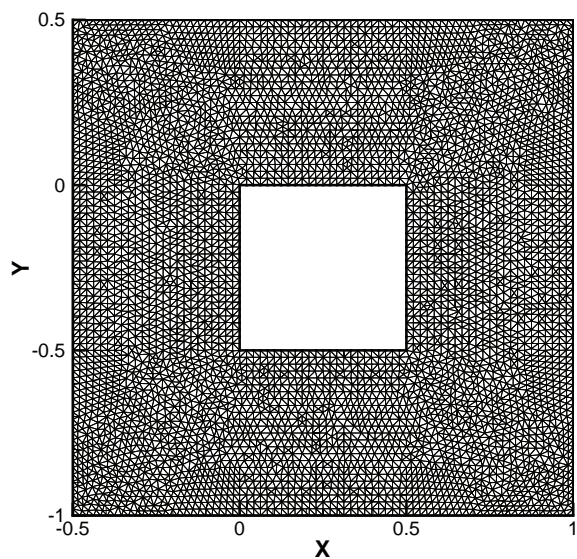}\\
	\caption{The mesh in the multi-connected domain.}\label{mesh}	
	\end{center}
\end{figure}

\begin{figure}[ht]
		\begin{center}
		\subfigure[$\kappa=1$]{\includegraphics[width=0.4\textwidth, height=56mm]{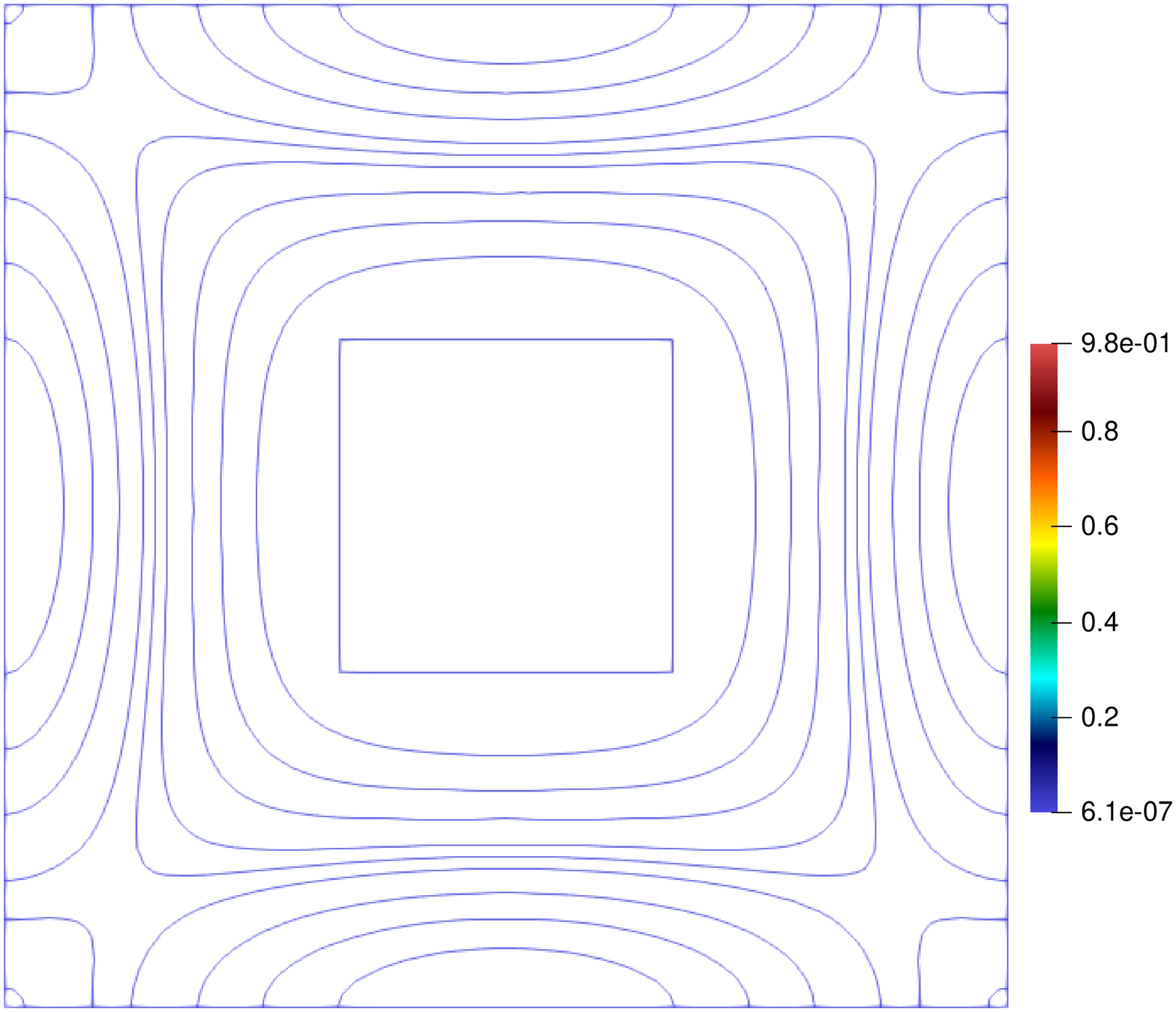}} 
		\subfigure[$\kappa=10$]{\includegraphics[width=0.4\textwidth, height=56mm]{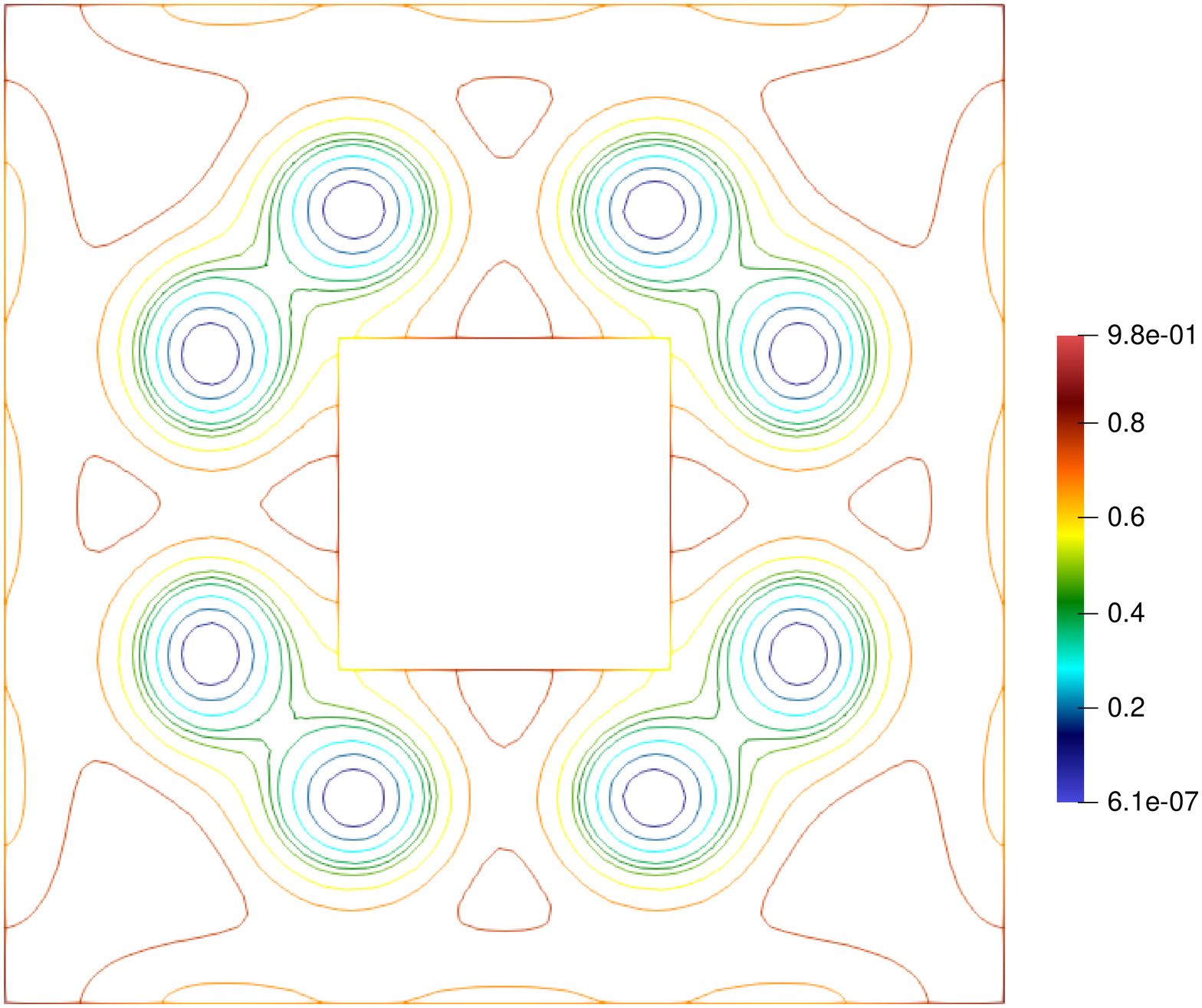}}\\
		 \subfigure[$\kappa=20$]{\includegraphics[width=0.4 \textwidth, height=56mm]{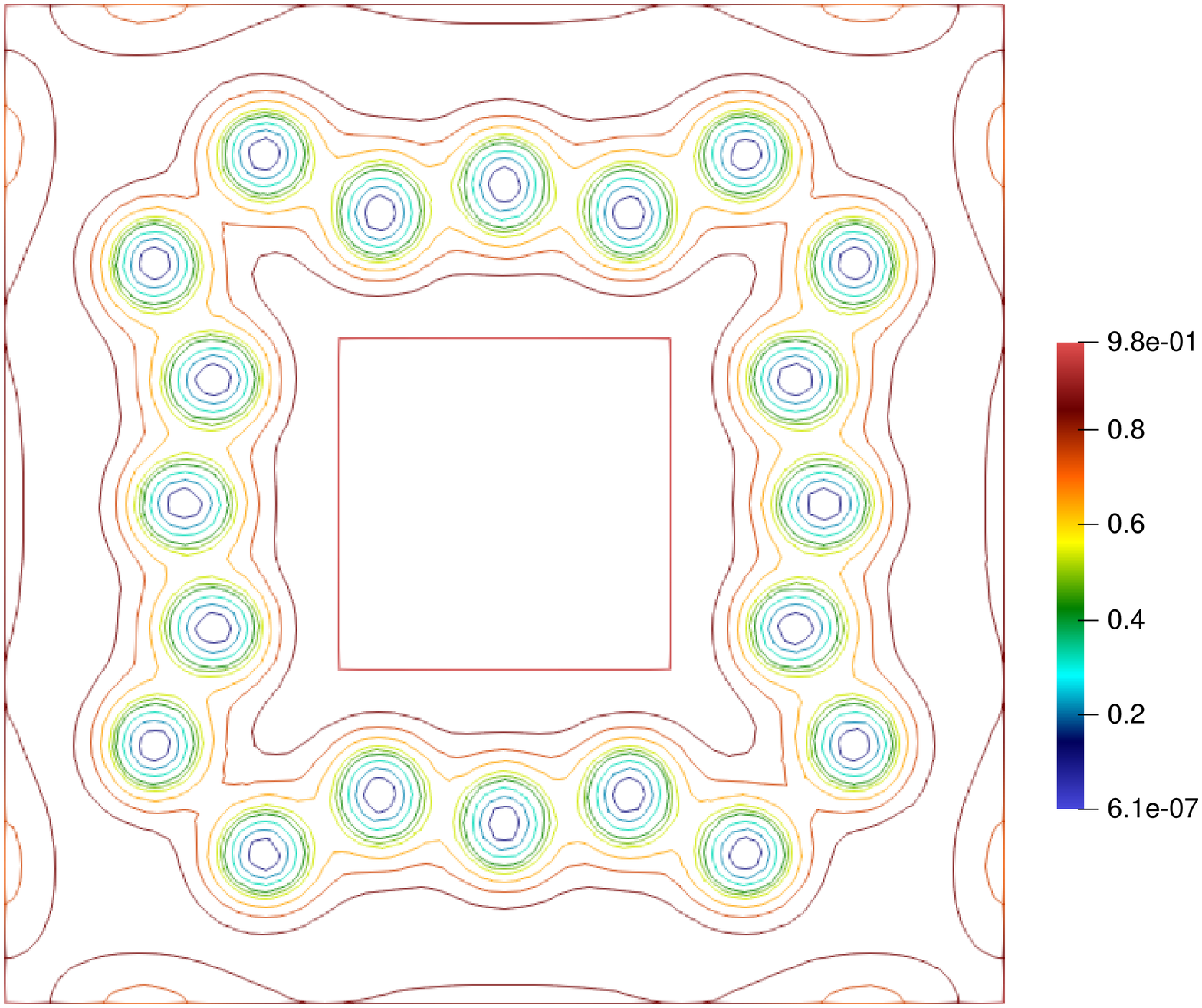}} 
		 \subfigure[$\kappa=30$]{\includegraphics[width=0.4 \textwidth, height=56mm]{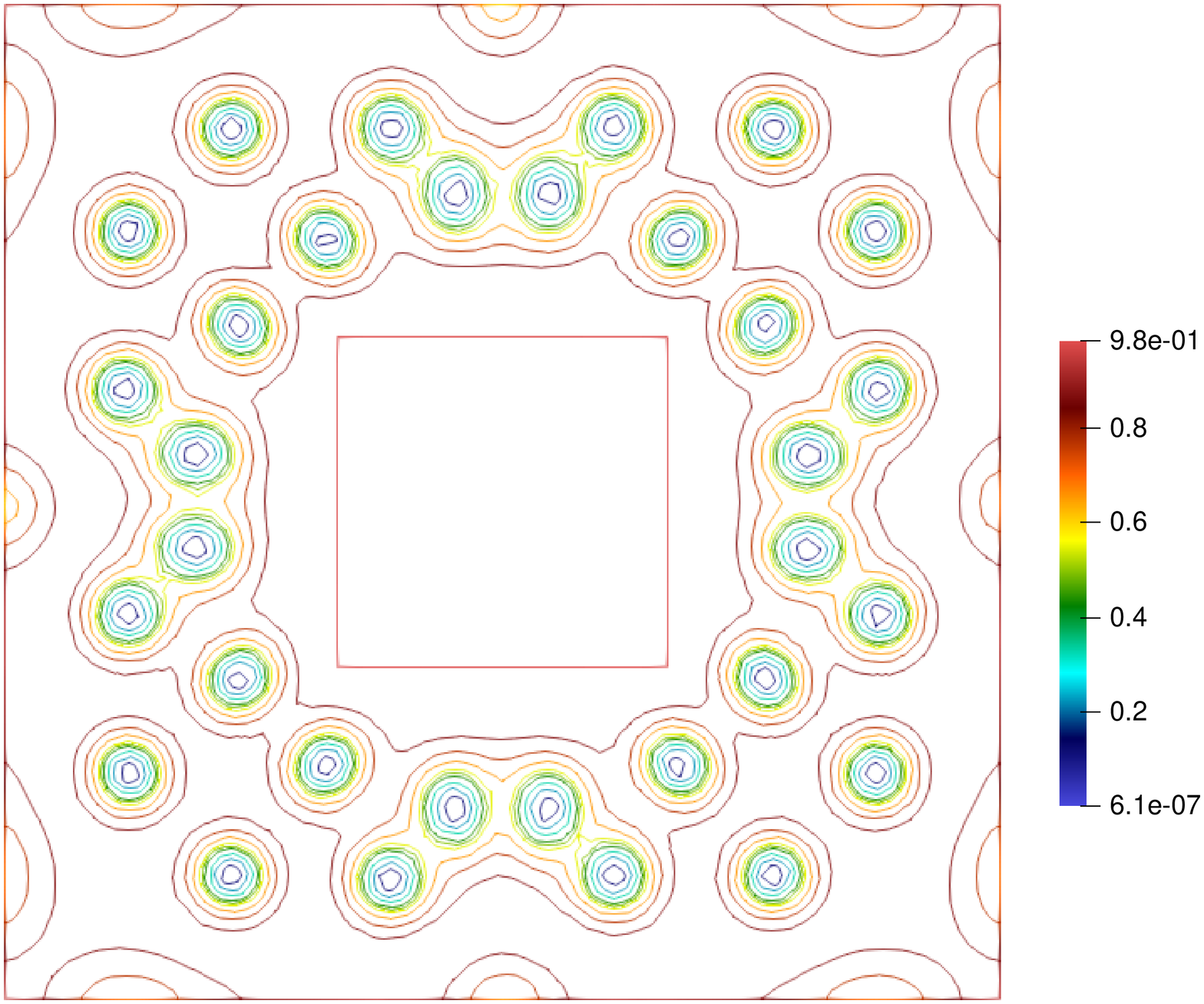}} 
		\caption{Contour plots of $|\psi|$ in the multi-connected domain with different $\kappa$.}\label{F2-2}
		\end{center}
	\end{figure}

\begin{figure}[ht]\begin{center}
		\includegraphics[width=0.66\textwidth, height=56mm]{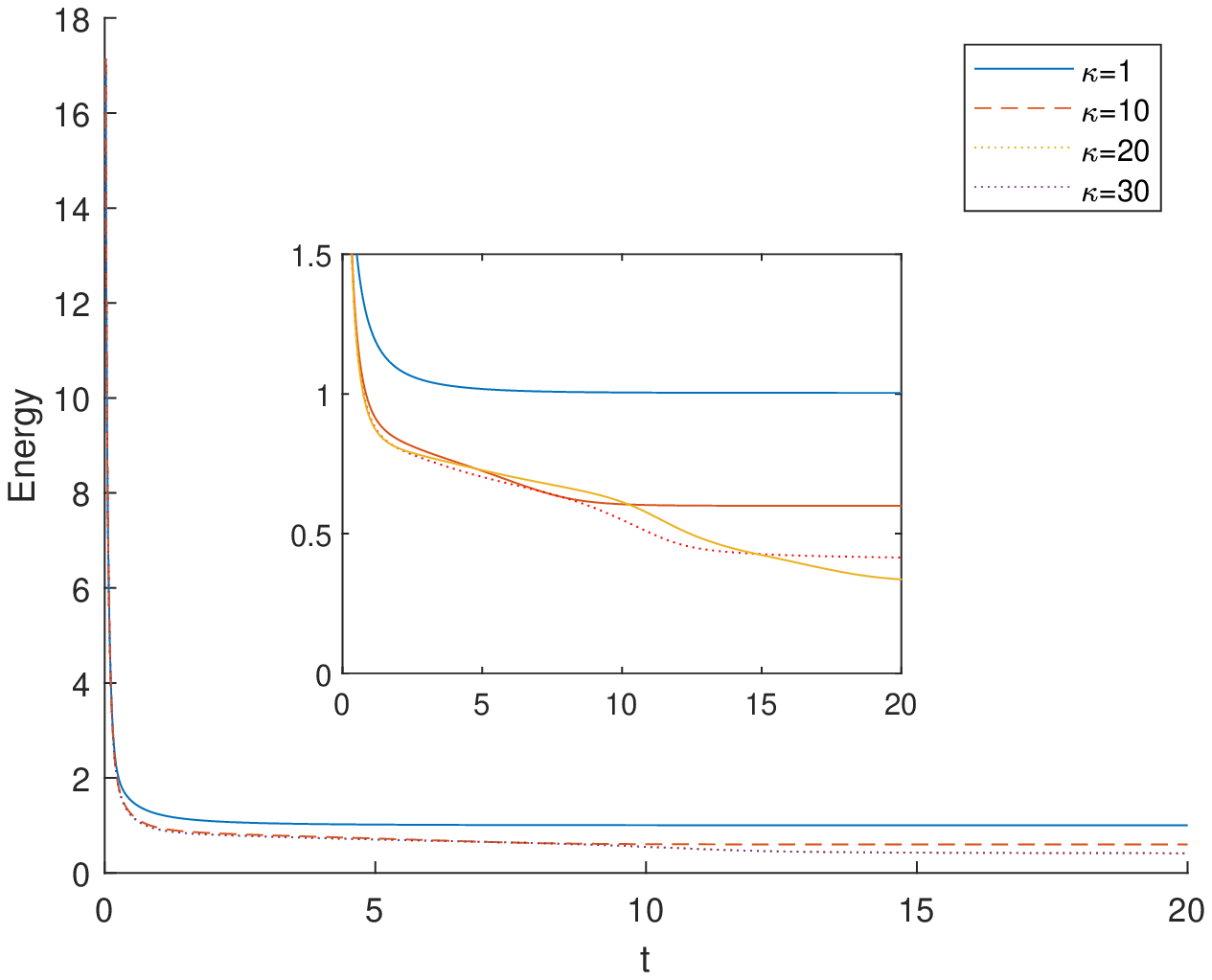}\\
		\caption{Time evolution of the discrete energy $\mathcal{G} (\psi,\A)$.}\label{F2-3}	
		\end{center}
	\end{figure}

\begin{figure}[ht]\begin{center}
		\includegraphics[width=0.66\textwidth, height=56mm]{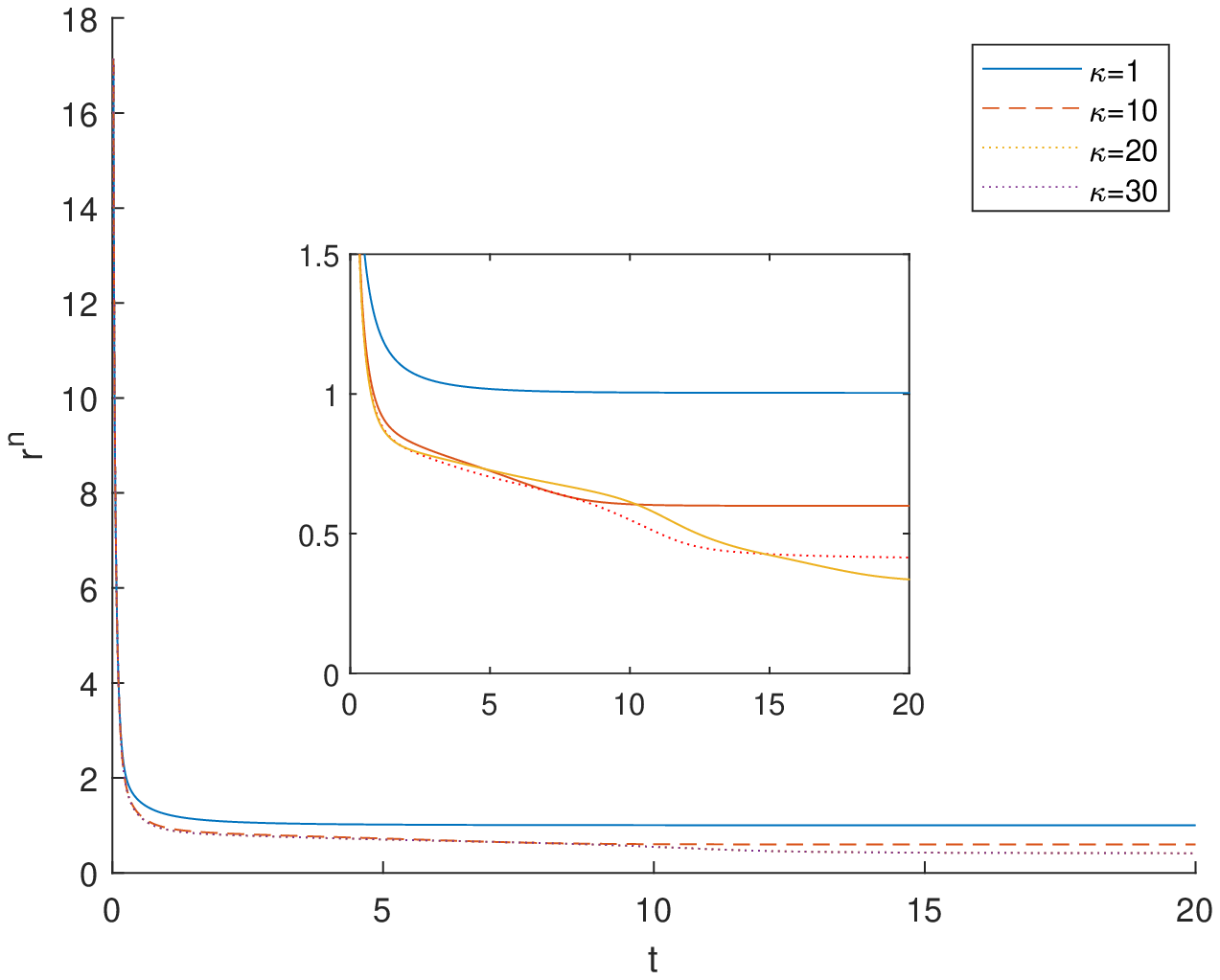}\\
		\caption{Time evolution of SAV $r$.}\label{F2-3b}	
		\end{center}
	\end{figure}

\begin{figure}[ht]\begin{center}
		\includegraphics[width=0.66\textwidth, height=56mm]{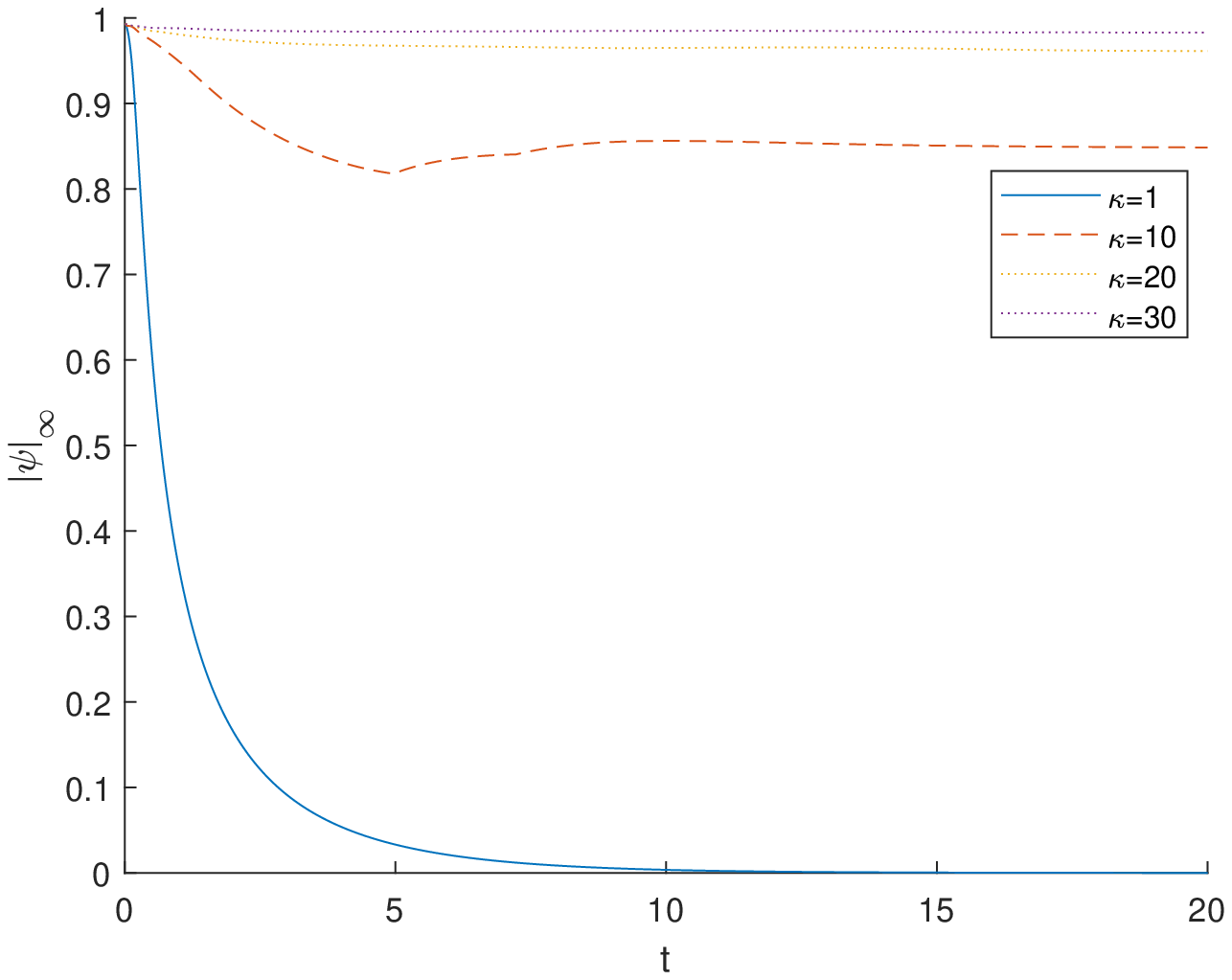}\\
		\caption{Time evolution of the $|\psi|_{\infty}$.}\label{F2-4}	
		\end{center}
	\end{figure}

This subsection shows some numerical results on the vortex simulation
 of the time-dependent GL equation in a multi-connected domain.
 We choose the domain $\Omega=\Omega_1/\Omega_2$, where $\Omega_1=[-0.5,1]\times[-1,0.5]$ and $\Omega_2=[0,0.5]\times[-0.5,0]$.
 We set the GL parameter $\kappa = 1, 10, 20$ and $30$, $\eta=1$. The initial conditions are set the same as in the first example. The mesh was given in Figure \ref{mesh}.
 Here, we choose the finite element as $\mathcal{V}_h^2, ~\mathbf{V}_h^2$ for $\psi$ and $\A$. The time step sizes are chosen as $\tau=0.01$ and the final time is $T=20$. The applied magnetic field $H = 5.0$. Figure \ref{F2-2} shows the contour plots of $|\psi|$ with different $\kappa$. It shows that the generalized SAV algorithm for the time-dependent GL equation is stable. Figure \ref{F2-3} presents the time evolution of the energy, we can see that the energy is reduced with the time $t$. Figure \ref{F2-3b} presents the time evolution of the SAV $r^n$, we can see that the SAV $r^n$ are reduced with the time $t$. Figure \ref{F2-4} shows the time evolution of maximum bound of magnetic parameter $|\psi|_{\infty}$. It shows that the maximum bound of the magnetic parameter is no bigger than $1$. It confirms the theoretical result.

	
\end{document}